\documentclass{amsart}

\usepackage{pdfpages}
\usepackage[procnames]{listings}
\usepackage{color}
\usepackage{cite}
\usepackage{relsize}
\usepackage{ amssymb }
\numberwithin{equation}{section}
\usepackage{amsmath}
\usepackage{mathtools}
\usepackage{amssymb}
\author{ Igor Rivin and Naser T. Sardari}
\title{Quantum Chaos on random Cayley graphs of ${\rm SL}_2[\mathbb{Z}/p\mathbb{Z}]$}
\date{\today}

	\newtheorem{thm}{Theorem}[section]

        \newtheorem{conj}[thm]{Conjecture}

	\newtheorem{rem}[thm]{Remark}	
	\newtheorem{lem}[thm]{Lemma}

	\newtheorem{cor}[thm]{Corollary}
	
	\theoremstyle{defi}

	\theoremstyle{pf}

	\numberwithin{equation}{section}

\begin {document}
\maketitle
\begin{abstract}
We  investigate the statistical behavior of the eigenvalues and diameter of  random Cayley graphs of  ${\rm SL}_2[\mathbb{Z}/p\mathbb{Z}]$ 
as the prime number $p$ goes to infinity. We  prove a density theorem for  the number of exceptional eigenvalues  of random Cayley graphs i.e. the eigenvalues with absolute value bigger than the optimal spectral bound.  Our numerical results suggest that random Cayley graphs of  ${\rm SL}_2[\mathbb{Z}/p\mathbb{Z}]$ and the explicit LPS Ramanujan  projective  graphs of $\mathbb{P}^1(\mathbb{Z}/p\mathbb{Z})$   have optimal spectral gap and diameter as the prime number $p$ goes to infinity. 
\end{abstract}
\tableofcontents
\section{Introduction}
\subsection{Motivation}
Our numerical experiments with random Cayley  graphs of ${\rm SL}_2[\mathbb{Z}/p\mathbb{Z}]$ are motivated by the experiments of  Lafferty and Rockmore \cite{LR3}, \cite{LR4},  \cite{LR} and later joint by Gamburd \cite{GLR} on the spectrum of the random Cayley graphs. In particular,  the authors investigated   the distribution of the normalized  consecutive distance (level spacing distribution) of eigenvalues of 4-regular random Cayley graphs of $SL_2(\mathbb{Z}/p\mathbb{Z})$,  symmetric group $S_n$ and large cyclic groups. In their initial papers,  they observed that the level spacing distribution is well modeled by Poisson distribution and hence do not follow  the distribution of the large symmetric random matrices (the Gaussian Orthogonal Ensemble (GOE)). This contradicts the general belief about the universality of the distribution of the level spacing of eigenvalues in the theory of  Quantum Chaos \cite{Berry} and  the strong numerical results for the same experiments with the 4-regular random graphs \cite{JMRR}, where the authors numerically showed  that the level spacing of the eigenvalues of a large 4-regular random graph is well modeled by the GOE distribution. Later, Lafferty and Rockmore in their  joint work with Gamburd \cite{GLR}  clarified this point by considering the eigenvalues associated to a single irreducible representation of ${\rm SL}_2[\mathbb{Z}/p\mathbb{Z}]$ and   using 6  instead of 4 generators for ${\rm SL}_2[\mathbb{Z}/p\mathbb{Z}]$; see \cite[Page 56]{Gamburd1999} for why 4 generators is not enough.  

We  partition  the spectrum of the random Cayley graph into the union of eigenvalues associated to irreducible representations of ${\rm SL}_2[\mathbb{Z}/p\mathbb{Z}]$.  This process is known as the desymmetrization of the spectrum that we explain in the next section.  We call the eigenvalues associated to a single irreducible representation of ${\rm SL}_2[\mathbb{Z}/p\mathbb{Z}]$, the monochromatic eigenvalues of the random Cayley graph. In other words,   Lafferty and Rockmore and Gamburd numerically showed that the level spacing of the monochromatic eigenvalues  of random Cayley graph of   ${\rm SL}_2[\mathbb{Z}/p\mathbb{Z}]$  approaches the GOE/GSE law of random matrix theory according to the parity of the irreducible representation of  ${\rm SL}_2[\mathbb{Z}/p\mathbb{Z}]$. Moreover, they observed that the level spacing distribution of the monochromatic eigenvalues of LPS Ramanujan graphs is well modeled by Poisson distribution.  This different statistical behavior  is explained by the arithmetic Quantum Chaos \cite{AQC}. In this paper, we  observe that the statistics of different monochromatic eigenvalues are essentially independent of each other and this is the reason that the level spacing of the eigenvalues of random Cayley graph without desymmetrization follows poisson distribution.        The main intuition behind the universality of  the distribution of the level spacing of the monochromatic eigenvalues comes from Quantum Chaos. 
%
Graphs are toy-models for  conjectures in Quantum Chaos where we can do numerical experiments on these finite combinatorial objects. 

In this paper, we investigate the universal properties of the monochromatic eigenvalues of large random Cayley graphs. The analogues results  are known in random matrix theory \cite{Deift06universalityfor}. More precisely,  we  investigate the statistical behavior of the eigenvalues of the fixed degree random Cayley graphs on  ${\rm SL}_2[\mathbb{Z}/p\mathbb{Z}]$ 
as prime number $p$ goes to infinity. We define a natural coloring of the eigenvalues of the adjacency matrix.  The colors are associated to the irreducible representations of the symmetry group of the graph that is  ${\rm SL}_2[\mathbb{Z}/p\mathbb{Z}]$. We define the monochromatic eigenvalues of the graph to be the eigenvalues associated to a single irreducible representation.  We prove the analogue of Wigner's semicircle law theorem \cite{Wigner} for monochromatic eigenvalues of ${\rm SL}_2[\mathbb{Z}/p\mathbb{Z}]$. Namely, we prove that the distribution of the monochromatic eigenvalues converges  to the Kesten-Mckay law.  Motivated by the universality of the asymptotic of the top eigenvalue 
of the random matrices, we investigate the asymptotic of the top monochromatic eigenvalues of the random $d$-regular Cayley graphs on ${\rm SL}_2[\mathbb{Z}/p\mathbb{Z}]$ when $d$ is fixed and the  prime number $p$ goes to infinity. 
Our numerical results suggest that the top nontrivial eigenvalue of the random $d$ regular Cayley graph  on ${\rm SL}_2[\mathbb{Z}/p\mathbb{Z}]$ is optimal and is bounded by the Ramanujan bound $2\sqrt{d-1}+\varepsilon$ almost surely as the prime number $p$ goes to infinity. For each family of monochromatic eigenvalues,    we show that the number of monochromatic eigenvalues that are bigger than $2\sqrt{d-1}+\varepsilon$ decay exponentially as $p$ goes to infinity.  Our result is uniform and does not depend on the color of the monochromatic eigenvalues.   We call this a density theorem for the number of exceptional monochromatic eigenvalues of random Cayley graph  on ${\rm SL}_2[\mathbb{Z}/p\mathbb{Z}]$. Based on our numerical results  and this density theorem, we conjecture that random Cayley graphs of ${\rm SL}_2[\mathbb{Z}/p\mathbb{Z}]$ are almost Ramanujan graphs. This generalizes the celebrated result of Bourgain and Gamburd on the expansion of random Cayley graphs on ${\rm SL}_2[\mathbb{Z}/p\mathbb{Z}]$ \cite{BG}.  As a consequence of the work of  Lubetzky and Peres \cite{Peres} or the second author \cite{Naser}, our conjecture implies that a random symmetric set of generators with $2d$ elements  represents  almost all elements of  ${\rm SL}_2[\mathbb{Z}/p\mathbb{Z}]$  with the optimal word length $(1+\varepsilon )\log_{2d-1}(|{\rm SL}_2[\mathbb{Z}/p\mathbb{Z}]|)$. In other words, the essential diameter of large random Cayley graphs is optimal where the essential diameter of a family of graphs defined on ${\rm SL}_2[\mathbb{Z}/p\mathbb{Z}]$ is  $h_p$ if the distance of almost all pairs of the vertices in ${\rm SL}_2[\mathbb{Z}/p\mathbb{Z}]$ is less than $h_p$ as $p\to \infty$. In fact,  we give strong numerical evidence that the diameter of the $d$-regular random Cayley graph is optimal and is asymptotic  to $\log_{d-1} (|G|)$ where $|G|$ is the number of vertices of graph $G$.  This conjecture and our numerical results clarify  Lafferty and Rockmore speculation about the  limiting distribution of the top eigenvalue of random Cayley graph of ${\rm SL}_2[\mathbb{Z}/p\mathbb{Z}]$ as $p$ goes to infinity. We discuss this further in the next paragraph.

Lafferty and Rockmore predicted that  random Cayley graphs of $SL_{2} (\mathbb{Z}/p\mathbb{Z})$ or Cayley graphs of $SL_{2} (\mathbb{Z}/p\mathbb{Z})$ with fixed generator set in $SL_{2} (\mathbb{Z})$  are not Ramanujan graphs and in fact they asked about the limiting distribution of the top eigenvalue as $p$ goes to infinity; see \cite{LR3}[Page 133]. Their prediction is based on their numerical experiments with prime numbers of small size $p \leq 500$. In section (\ref{top}), we numerically investigate the asymptotic  of the top nontrivial monochromatic eigenvalues associated to the principal series  representations of  $SL_{2} (\mathbb{Z}/p\mathbb{Z})$ of both families of random and fixed generator Cayley graphs of $SL_{2} (\mathbb{Z}/p\mathbb{Z})$  where  $p\approx 4409$. We give numerical evidence that the families of  fixed or random  Cayley graphs of  ${\rm SL}_2[\mathbb{Z}/p\mathbb{Z}]$   are almost Ramanujan graphs  for large prime numbers $p$. In other words, the top eigenvalue of the family of $d$-regular random or fixed generator set of Cayley graphs of  ${\rm SL}_2[\mathbb{Z}/p\mathbb{Z}]$ converges to the Ramanujan bound $2\sqrt{d-1}$ as $p$ goes to infinity. 

More generally, let $X:={\rm SL}_2[\mathbb{Z}/p\mathbb{Z}]/H$ where $H$ is a subgroup of ${\rm SL}_2[\mathbb{Z}/p\mathbb{Z}]$. We note that ${\rm SL}_2[\mathbb{Z}/p\mathbb{Z}]$ acts from the left on $X$. Let $S$ denote a finite subset of ${\rm SL}_2[\mathbb{Z}/p\mathbb{Z}]$. The Schreier graph defined on $X$ associated to the generating set $S$ is the graph with vertex set $X$ such that $x\in X$ is connected to $s.x$ for all $s\in S$.  We note that the Schreier graph associated to the pair $(X,S)$ is the quotient of the Cayley graph $({\rm SL}_2[\mathbb{Z}/p\mathbb{Z}],S)$  by $H$. Therefore,  the spectrum of Schreier graphs are a subset of the associated Cayley graph. As a result, every random or fixed generator Schreier graph associated to  ${\rm SL}_2[\mathbb{Z}/p\mathbb{Z}]$ is almost Ramanujan graph for large prime numbers. We also study the fluctuation of the top monochromatic eigenvalues around the Ramanujan bound $2\sqrt{d-1}$  for random and fixed generator Cayley graphs on   ${\rm SL}_2[\mathbb{Z}/p\mathbb{Z}]$.    The top nontrivial eigenvalue of a random Cayley graph on $SL_{2} (\mathbb{Z}/p\mathbb{Z})$ is the maximum over all top monochromatic eigenvalues. 
 Since the number of the monochromatic eigenvalues of $SL_{2} (\mathbb{Z}/p\mathbb{Z})$ is linear in $p$ and heuristically different color eigenvalues are independent of each other then with probability 1 one of them is bigger than the Ramanujan bound $2\sqrt{d-1}$  as $p\to \infty$. On the other hand, the random projective Cayley graphs $P^1(\mathbb{Z}/p\mathbb{Z})$ are Ramanujan graphs with a positive probability. Since, the nontrivial eigenvalues of the projective Cayley graph and monochromatic and are associated to  Steinberg representation of  $SL_{2} (\mathbb{Z}/p\mathbb{Z})$. In section (\ref{top}), we verify numerically the last two observation on the top non-trivial eigenvalues of Cayley graph of $SL_{2} (\mathbb{Z}/p\mathbb{Z})$ and the projective Cayley graph of $P^1(\mathbb{Z}/p\mathbb{Z})$.

We conduct similar experiments with  the random $2d$-regular  Schreier graphs of the Symmetric group $S_n$  and Schreier graph of ${\rm SL}_3[\mathbb{Z}/p\mathbb{Z}]$. Our results show that the random Cayley graphs on these groups have optimal spectral gap. Note that the probability space of the random $2d$-regular  Schreier graphs of the Symmetric group $S_n$ on $\{1,2,\dots,n\}$  is contiguous with the uniform distribution on the random $2d$-regular graphs with $n$ vertices \cite[Page 497]{HLW}. Two sequence of measures are  contiguous if asymptotically they share the same support. Therefore, almost Ramanujan property holds with probability 1 under both probability measures.  Our numerical results in this special family of Schreier graphs on $S_n$ is  reduced to the numerical results of  Jakobson, Miller, Rivin, and Rudnick \cite{JMRR} on the level spacing distribution of the random $d$ regular graphs and the numerical results of Miller, Novikoff and Sabelli~\cite{miller2008} on the  distribution of the top nontrivial eigenvalue of the random  $d$ regular graphs. The optimality of the top eigenvalue and the diameter of these random Schreier graphs  are reduced to Friedman's second eigenvalues theorem and Bollobas-Vega's theorem for the diameter of random regular graphs.

\subsection{Desymmetrization of the spectrum}

\noindent 
 In this section, we introduce a natural coloring of the spectrum of the adjacency matrix  of a Cayley graph.  The colors are associated to the irreducible representations of the associated group. We define the monochromatic eigenvalues of the graph to be the eigenvalues associated to a single irreducible representation of the group. We denote this process by the desymmetrization of the spectrum of the Cayley graph.

  Let $G$ be a finite group and $S:=\{g_1^{\pm}, \dots, g_d^{\pm}   \}$ a symmetric generating  set of $G$ with $2d$ elements. Let $X$ be the $2d$-regular Cayley graph of $G$ generated by $S$. Let $\mathbb{C}(X)$ be the space of complex valued functions on the graph $X$. Let $\Delta$ be the Laplacian operator defined on $\mathbb{C}(X)$ by 
 $$\Delta (f) (x)= \sum_{y\sim x} f(y),$$
 where $f\in \mathbb{C}(X) $ and $y\sim x$ if there is an edge between $y$ and $x$ in $X$. We identify $\mathbb{C}(X)$ by the complex group ring $\mathbb{C}(G)$. We note that the Laplacian operator on $\mathbb{C}(X)$ is associated to right multiplication by $$g_1+g_1^{-1}+\dots+g_d+g_d^{-1}\in \mathbb{C}(G). $$
 As a result the eigenvalues of the Laplacian operator $\Delta$ on  $\mathbb{C}(X)$ is associated the eigenvalues of the right multiplication by $g_1+g_1^{-1}+\dots+g_d+g_d^{-1}\in \mathbb{C}(G). $  It is a classical theorem that the regular representation of the finite group $G$ decomposes into a sum over all irreducible representations $\pi$ of $G$ with multiplicity $d(\pi)$ where $d(\pi)$ is the dimension of the irreducible representation $\pi$.  Hence the eigenvalues of the Laplacian operator $\Delta$ are the union over all eigenvalues of the linear operator
 $$\pi(g_1)+\pi(g_1^{-1})+\dots+ \pi(g_d)+\pi(g_d^{-1}), $$
 with multiplicity $d(\pi)$.
 %
 Let $\{\lambda_1, \dots, \lambda_n\}$ denote the spectrum of Laplacian operator  $\Delta$. We color the eigenvalue $\lambda_{i}$ with the associated  irreducible representation $\pi$ of $G$. We denote this coloring the desymmetrization of the spectrum.  We define the monochromatic eigenvalues of the graph $X$ to be the subset of eigenvalues associated to a single irreducible representation $\pi$.   Note that if $X$ is a Schreier graph on a finite group $G$, then the spectrum of $X$ are a subset of the spectrum of the associated Cayley graph on $G$. As a result, this embedding gives a natural coloring of the spectrum of the Schreier graph $X$ by the irreducible representations of $G$. 
 
 In this paper, we work with Schreier graphs on ${\rm SL}_2[\mathbb{Z}/p\mathbb{Z}]$. By desymmetrization of the spectrum, we mean the coloring of the spectrum given by irreducible representations of ${\rm SL}_2[\mathbb{Z}/p\mathbb{Z}]$.    Our experiments  with the  families of Schreier graphs show that different sets of monochromatic eigenvalues are statistically  independent of each other. On the other hand, for the generic families of Schreier graphs the level spacing distribution of the monochromatic eigenvalues  
are well-modeled by  the law of random matrix theory (e.g. GOE/GSE). In particular they repel each other and behave differently from the statistics of the random sample of independent points (Poisson process). This is consistent with the prediction of Quantum Chaos for the distribution of eigenvalues  in the high energy limit. 
\subsection{Statement of results and outline of the paper}
In this section, we introduce the families of Schreier graphs of ${\rm SL}_2[\mathbb{Z}/p\mathbb{Z}]$  that we investigate in this paper. We observe two different statistical behaviors of the monochromatic eigenvalues that is explained by the arithmetic and  generic Quantum Chaos. The arithmetic behavior is observed for  the family of  LPS Ramanujan graphs.  The generic behavior is observed  in  different natural families of Schreier graphs. It includes  random Schreier graphs given by the action of the Symmetric group $S_n$ on $\{1,2,\dots,n\}$, the  random Schreier graph on ${\rm SL}_2[\mathbb{Z}/p\mathbb{Z}]/H$, where $H$ is a fixed algebraic subgroup of ${\rm SL}_2[\mathbb{Z}/p\mathbb{Z}]$ and the generators are chosen randomly according to the uniform measure   on  ${\rm SL}_2[\mathbb{Z}/p\mathbb{Z}]$, the Schreier graph on ${\rm SL}_2[\mathbb{Z}/p\mathbb{Z}]$ with respect to  fixed generic elements in ${\rm SL}_2[\mathbb{Z}]$ as $p$ goes to infinity or even the Schreier graph ${\rm SL}_2[\mathbb{Z}]/H$  with respect to the LPS generators, where $H$ chosen uniformly among subgroups of index $ n$ in ${\rm SL}_2[\mathbb{Z}]$ as $n$ goes to infinity. 

 Next, we introduce the families of the  Schreier graphs on the left cosets of  ${\rm SL}_2[\mathbb{Z}/p\mathbb{Z}]/H$, where $H$ is a subgroups of ${\rm SL}_2[\mathbb{Z}/p\mathbb{Z}]$. If $H$ is the subgroup of the upper triangular matrices $B$ (Borel subgroup) of ${\rm SL}_2[\mathbb{Z}/p\mathbb{Z}]$, then we call the Schreier graph by the projective  Cayley graph. Note that the action of ${\rm SL}_2[\mathbb{Z}/p\mathbb{Z}]$ from the left on ${\rm SL}_2[\mathbb{Z}/p\mathbb{Z}]/B$ is isomorphic to the action of ${\rm SL}_2[\mathbb{Z}/p\mathbb{Z}]$ on $\mathbb{P}^1(\mathbb{Z}/p\mathbb{Z})$ by  Mobius transformation. The eigenvalues are colored by the trivial and the Steinberg representation of  ${\rm SL}_2[\mathbb{Z}/p\mathbb{Z}]$. The nontrivial eigenvalues of the projective   Schreier graph are monochromatic.  If $H$ is the subgroup of unipotent matrices $N$ of ${\rm SL}_2[\mathbb{Z}/p\mathbb{Z}]$, then we denote this Schreier graph by  linear Cayley graph. The action of  ${\rm SL}_2[\mathbb{Z}/p\mathbb{Z}]$ from the left on ${\rm SL}_2[\mathbb{Z}/p\mathbb{Z}]/N$ is isomorphic to the action of ${\rm SL}_2[\mathbb{Z}/p\mathbb{Z}]$ by linear maps on $\mathbb{A}^2(\mathbb{Z}/p\mathbb{Z})-\{(0,0) \}$. The monochromatic eigenvalues of  linear Schreier graph are associated to principal series representations of ${\rm SL}_2[\mathbb{Z}/p\mathbb{Z}]$. If the generators are uniformly chosen   from  ${\rm SL}_2[\mathbb{Z}/p\mathbb{Z}]$, we denote this by random linear Cayley graph. If the generators are LPS generators then we denote the graph by LPS Schreier graph. If the generators are fixed inside ${\rm SL}_2[\mathbb{Z}]$ we call the Schreier graphs the fixed generators Schreier graphs. Finally, we denote the  Schreier graphs of the Symmetric group $S_n$ on $\{1,2,\dots,n\}$ with respect to symmetric set of generators by the permutation graphs. 

 In section \ref{dist}, we prove that  the bulk  distribution of the monochromatic eigenvalues of random Cayley graphs converges weakly to Kesten-Mckay law. More precisely in Theorem~(\ref{ddd}),  we show that the discrepancy of the distribution of the monochromatic eigenvalues associated to $\pi$ and Kesten-Mckay law is $O(1/\log(p))$ for almost all  generators as $p$ goes to infinity.  We denote the discrepancy between the measures $\nu$ and $\mu$ by  $D(\nu ,\mu)$, where 
 $$D(\nu,\mu) = \sup\{ |\nu(I)-\mu(I)|: I =[a,b]\subset \mathbb{R}\}.$$
 According to the random matrix prediction and our numerical results, the discrepancy of generic element should be $O(\log(d(\pi))/d(\pi))$ where $d(\pi)=O(p)$ is the dimension of the irreducible representation $\pi$. 
 \begin{thm} \label{ddd} \normalfont  Fix $d\geq 2$. Let $\pi\neq 1$ be any nontrivial  irreducible representation of $G=SL_2(\mathbb{Z}/p\mathbb{Z})$. Let $S:=\{ g_1^{\pm},\dots,g_d^{\pm}\}$ be a random symmetric subset of $SL_2(\mathbb{Z}/p\mathbb{Z})$ chosen uniformly. Let   $\mu_{\pi,S}$ be the spectral probability measure of $ \pi (g_1+g_1^{-1}+ \dots + g_d+g_d^{-1})$ and $f_{2d}$ be the Kesten-Mckay law with parameter $2d$. Let $D(\mu_{\pi,S}, f_{2d})$ be the discrepancy between $\mu_{\pi,S}$ and $f_{2d}$.  Then 
$$D(\mu_{\pi,S}, f_{2d}) \ll \frac{1}{\log (p)},$$
 for a symmetric subset $S\subset G $  with probability $(1-O(p^{-1/4}))$. The constants involved in $\ll $ and $O(p^{-1/4})$ only depends on $d$. 
\end{thm}

 We say  an eigenvalue $\lambda$ of $\pi (g_1+g^{-1}_1+\dots+ g_{d}+g^{-1}_d)$ is an exceptional eigenvalue if $\lambda >  |2\sqrt{2d-1}| $. Let $N(\alpha, \pi,S)$ be the number of monochromatic exceptional eigenvalues such that $|\lambda| > \alpha (2\sqrt{2d-1})$ where $\alpha >1$. As a corollary of the method of the proof of Theorem~(\ref{ddd}),  we prove a density theorem for the number of  exceptional monochromatic eigenvalues of random Cayley graphs.
 Note that the number of monochromatic eigenvalues is $d(\pi)= O(p)$. So, the trivial upper bound is $N(\alpha, \pi,S) =O(p)$.    In the following corollary, we give a sub-exponential  upper bound on $N(\alpha, \pi,S)$. 
\begin{cor}\label{cccc}\normalfont \noindent 
Let $N(\alpha, \pi,S)$ be the number of exceptional monochromatic eigenvalues associated to $\pi$ such that $|\lambda| > \alpha (2\sqrt{2d-1})$ where $\alpha >1$. Then
$$N(\alpha, \pi,S) \ll p^{1-\frac{(\log_{2d}{\alpha})+o(1)}{2} }. $$
\end{cor}

This corollary is  the analogue of the Iwaniec's density theorem for the number of exceptional eigenvalues of the laplacian for congruence subgroups; see \cite{IS}. By Selberg's $1/4$ conjecture there is no exceptional eigenvalue for congruence subgroups \cite{Selberg}.    In section~\ref{top}, we present our numerical results that show  the nontrivial eigenvalues are bounded by 
$$|\lambda| < 2(1+\varepsilon)\sqrt{d-1},$$ 
for almost all symmetric generator sets. Our numerical results together with this density theorem suggest that random Cayley graphs of $SL_2(\mathbb{Z}/p\mathbb{Z})$ are almost Ramanujan. This generalizes the celebrated result of Bourgain and Gamburd on the expansion of random Cayley graphs.  Consequently,  the essential  diameter of  random Schreier graphs on ${\rm SL}_2[\mathbb{Z}/p\mathbb{Z}]$  is optimal as $p$ goes to infinity. The essential diameter of a graph is $d$ if $99\%$ of the distance of pairs of vertices is less than $d$.

In section \ref{diam},  we present our numerical experiments on the diameter of the fixed and random generator families of Cayley and projective Cayley  graphs. Our experiments  show that random generator families of Cayley and projective Cayley graphs have optimal diameter. On the other hand, the diameter of the Cayley graph with LPS generators that is a special fixed generator family of Cayley graphs is not optimal. The second named author in  \cite{Naser}, showed that the diameter of a  family  of LPS Ramanujan graphs is greater than $4/3\log_{d-1}(n)$  where $n$ is the number of vertices and $d$ is the degree of the graph. This lower bound is related to the repulsion of integral points lying on quadrics. The second named author formulated a conjecture on the optimal exponent for strong approximation for quadratic forms in 4 variables  \cite[Conjecture 1.9]{Optimal} that implies the diameter of LPS Ramanujan graphs is asymptotic to $4/3\log_{d-1}(n)$. More recently,  he developed and implemented a polynomial time algorithm for navigating LPS Ramanujan graphs under a polynomial time algorithm for factoring integers and an arithmetic conjecture on the distribution of numbers representable as a sum of two squares \cite{Sss}. The numerical results from that algorithm strongly support that the diameter of LPS Ramanujan graphs is asymptotic to $4/3\log_{d-1}(n)$; see \cite{Sss}.   This shows that the diameter of Random Cayley graphs  is asymptotically smaller than the diameter of the LPS Cayley graphs. On other hand, we present new numerical results for the diameter of the LPS projective Cayley graphs that show their diameter is optimal and there is no repulsion in the LPS Scherier graph defined on $\mathbb{P}^1(\mathbb{Z}/q\mathbb{Z})$. This is the first family of explicit graphs that have both optimal diameter and spectral gap.

%

In section~\ref{extend}, we generalize Theorem~{\ref{ddd}} to ${\rm SL}_n[\mathbb{Z}/p\mathbb{Z}]$ for fixed $n\geq 3$ as $p\to \infty$. We prove that the distribution of the monochromatic eigenvalues of $d$-regular random Cayley graphs of $G=SL_n(\mathbb{Z}/p\mathbb{Z})$ for any $n\geq 3$  converges to  Kesten-Mackay law for almost all random generators as $p$ goes to infinity. 

\begin{thm}\label{lasttheorem}  \normalfont  Fix $d\geq 2$. Let $\pi$ be any nontrivial  irreducible representation of $G=SL_n(\mathbb{Z}/p\mathbb{Z})$. Let $S:=\{ g_1^{\pm},\dots,g_d^{\pm}\}$ be a random symmetric subset of $G$ with respect to the uniform measure on finite group $G$. Let   $\mu_{\pi,S}$ be the spectral probability measure of $ \pi (g_1+g_1^{-1}+ \dots + g_d+g_d^{-1})$ and $f_{2d}$ be the Kesten-Mckay law with parameter $2d$ . Let $D(\mu_{\pi,S}, f_{2d})$ be the discrepancy between $\mu_{\pi,S}$ and $f_{2d}$.  Then 
$$D(\mu_{\pi,S}, f_{2d}) \leq \varepsilon,$$
 for symmetric subset $S\subset G $  with probability $(1-O_{\varepsilon}(p^{-1}))$. The constants in $O_{\varepsilon}(p^{-1})$ depends only on $n$, $d$  and $\varepsilon$.
\end{thm}

We formulate two conjectures that predict  the statistics of the monochromatic eigenvalues of Cayley graphs on  semi-simple groups defined over $\mathbb{Z}$ and the symmetric group $S_n$. Our conjecture for the semi-simple groups generalize the super-strong approximation. It states that in the family of fixed  or  random generator Cayley graphs on  semi-simple groups the top eigenvalue  is asymptotic to the Ramanujan bound as $p$ goes to infinity. 
%
%
%

\begin{conj}\label{ccon}
\normalfont  Let $G$ be a semi-simple algebraic group defined over $\mathbb{Z}$.
Fix $k\geq 3$. Let $p$ be a prime number and $\pi$ be any nontrivial irreducible representation of $G(\mathbb{Z}/p\mathbb{Z})$. Let  $\{g_1, \dots,g_k\}$ be a random subset of $G(\mathbb{Z}/p\mathbb{Z})$ chosen uniformly. The distribution of the eigenvalues of the  self-adjoint  operator $ \pi (g_1+g_1^{-1}+ \dots + g_k+g_k^{-1})$ converges weakly to the Kesten-Mckay law for almost all $\{g_1, \dots,g_k\}$   as $p$ goes to infinity. The level spacing distribution  of the monochromatic  eigenvalues  is well modeled by the the law of random matrix theory  and the top nontrivial eigenvalues is bounded by Ramanujan bound $(2+\varepsilon)\sqrt{2k-1}$ almost surly as $p$ goes to infinity.   Finally, the statistics of two different monochromatic eigenvalues of $g_1+g_1^{-1}+ \dots + g_k+g_k^{-1}$ associated to two irreducible represntations of $\pi_1$ and $\pi_2$ of $G(\mathbb{Z}/p\mathbb{Z})$ are independent of each other.  
\end{conj}
\begin{rem}
\normalfont If $G$ is not a semi-simple group then the top nontrivial eigenvalue is not asymptotically bounded by the Ramanujan bound $(2+\varepsilon)\sqrt{2k-1}$. In fact Klawe \cite{Maria}, showed that the random Cayley graph of the affine linear maps do not have spectral gap. The top nontrivial eigenvalues is asymptotic to  $2k$ as $p$ goes to infinity. On the other hand, from the Super-strong approximation theorem we have expansion property for random Cayley graphs of  $G(\mathbb{Z}/p\mathbb{Z})$. The condition $k\geq 3$ is necessary. Because  the centralizer subgroup of two generic element $g_1$ and $g_2$ is nontrivial and that changes the statistics of the monochromatic eigenvalues of $( g_1+g_1^{-1}+ g_2+g_2^{-1})$; see \cite[Page 56]{Gamburd1999} for why 4 generators is not enough.
\end{rem}

\begin{conj} \normalfont
 Let $S_n$ be the symmetric group. Fix $k\geq 2$. Let $\pi\neq \text{sign}$ be a nontrivial irreducible representation of $S_n$. Let  $\{g_1, \dots,g_k\}$ be a random subset of $S_n$ chosen uniformly. The distribution of  the eigenvalues of the self-adjoint  operator $\pi (g_1+g_1^{-1}+ \dots + g_k+g_k^{-1})$ converges weakly to the Kesten-Mckay law for almost all $\{g_1, \dots,g_k\}$  as $n$ goes to infinity. The level spacing of the eigenvalues is well modeled by GOE level spacing distribution and the top nontrivial eigenvalues is bounded by Ramanujan bound $(2+\varepsilon)\sqrt{2k-1}$ as $n$ goes to infinity. Finally, the statistics of two different monochromatic eigenvalues of $g_1+g_1^{-1}+ \dots + g_k+g_k^{-1}$ associated to two irreducible represntations of $\pi_1$ and $\pi_2$ of $S_n$ are independent of each other.    \end{conj}

\begin{rem}\normalfont
Let $\pi$ be the standard irreducible representation of $S_n$. Then the optimality of the top eigenvalue is the Friedman's second eigenvalues theorem. So, this conjecture extends Friedman's second eigenvalues theorem to all the nontrivial irreducible representation $\pi\neq \text{sign}$ of $S_n$.
\end{rem}

\section{The bulk  distribution of the monochromatic  eigenvalues }\label{dist}
In this section,  we give a proof of Theorem~(\ref{ddd}).
    Let $S:=\{ g_1^{\pm},\dots,g_d^{\pm}\}$ be a random symmetric subset of $SL_2(\mathbb{Z}/p\mathbb{Z})$ and   $G_{S}$ be Cayley graph of $SL_2(\mathbb{Z}/p\mathbb{Z})$ generated by $S$. As pointed out in the introduction, the logarithmic lower bound on the girth of  random Cayley graphs on $SL_{2} (\mathbb{Z}/p\mathbb{Z})$ \cite{GHSS} implies that the distribution of the eigenvalues of  $G_{S}$ converges to Kesten-Mckay law with parameter $2d$. We prove a refined version of this result. We recall that the Kesten-McKay's law with parameter $2d$  is  the spectral measure of the infinite $2d$-regular tree which has the following density:
\begin{equation*}
f_{2d}(x)= \frac{2d\sqrt{4(2d-1)-x^2}}{2\pi ((2d)^2-x^2)},   \text {     for    }  -2\sqrt{2d-1} \leq x \leq 2 \sqrt{2d-1}.
\end{equation*}
    Let  $\pi$ be an irreducible representation of $SL_{2} (\mathbb{Z}/p\mathbb{Z})$.  The monochromatic eigenvalues of Cayley graph $G_{S}$ associated to $\pi$ are the eigenvalues of  the  self-adjoint  operator $ \pi (g_1+g_1^{-1}+ \dots + g_d+g_d^{-1})$. We denote the bulk distribution  of monochromatic eigenvalues associated to $\pi$ and generator set $S$ by  $\mu_{\pi,S}$. Next, we give the proof of Theorem~(\ref{ddd}).

%
\begin{proof}[Proof of Theorem~(\ref{ddd})] \normalfont   
We begin by bounding the discrepancy  $D(\mu_{\pi,S}, f_{2d})$ with $|\mu_{\pi,S}(x^m)- f_d(x^m)|$ for $m\ll \log(p)$ where $\mu_{\pi,S}(x^m)$ and $f_{2d}(x^m)$   are the  the $m$th moment of $\mu_{\pi,S}$ and $f_{2d}$, respectively. This passage is straightforward by approximating the indicator function of an interval by Selberge polynomials; see\cite{vaaler1985}, \cite{MR3308963}.    We refer the reader to \cite[Proof of Theorem 1.3]{Gamburd1999}. In particular, if $S$ is a generator set for $SL_2[\mathbb{Z}/p\mathbb{Z}]$ such that for every integer  $m \leq \alpha \log (p)$, 
 \begin{equation}\label{mme}
|\mu_{\pi,S}(x^m)- f_d(x^m)|\ll p^{-\beta},
\end{equation} 
 for some $0<\alpha$ and  $0< \beta$.  Then 
 \begin{equation}\label{discp}D(\mu_{\pi,S}, f_{2d}) \ll \frac{1}{\log (p)}.\end{equation}

Therefore, it is enough to prove that for some $\alpha$ and $\beta$  inequality $(\ref{mme})$ holds for generator set $S$ with probability $1-O(p^{-1/4})$. Our strategy to prove inequality (\ref{mme}) is to exclude a measure $O(p^{-1/4})$ family of symmetric subsets $S$ and give an upper bound on the variance of the $m$th moment of the remaining symmetric sets $S$. Then by Chebyshev's inequality, we deduce that for some $\alpha$ and $\beta$ the   inequality (\ref{mme}) holds for  symmetric subsets $S$    with probability $(1-O(p^{-1/4}))$.  
\\
\\
First, we describe the subsets $S$  that we exclude and prove their probability is $O(p^{-1/4})$.  Let $\chi_{\pi}$ be the character associated to the irreducible representation $\pi$ of  $SL_{2} (\mathbb{Z}/p\mathbb{Z})$. The  $m$th moment of $\mu_{\pi,S}$ is
$$\mu_{\pi,S}(x^{m})=\frac{1}{d(\pi)}\chi_{\pi}((g_1+g_1^{-1}+ \dots + g_d+g_d^{-1})^{m}),$$
where $d(\pi)$ is the dimension of the irreducible representation $\pi$. We expand the product and write 
\begin{equation}\label{summ}\mu_{\pi,S}(x^{m})=\frac{1}{d(\pi)}\mathlarger{\sum}_{i_1,\dots,i_m}\chi_{\pi}( g_{i_1}\dots g_{i_m}),\end{equation}
where the sum is over all words of length $m$ with letters in $S:=\{ g_1^{\pm},\dots,g_d^{\pm}\}$.
 We denote $g_{i_1}\dots g_{i_m}$  by trivial word if $g_{i_1}\dots g_{i_m}=id$ as an element in the free group generated by formal letters $\{ g_1^{\pm},\dots,g_d^{\pm}\}$ and by nontrivial word otherwise. We split the sum~(\ref{summ}) into two parts over the trivial and nontrivial words. Let $N(d,m)$ be the number of trivial words of length $m$  with formal letters $\{ g_1^{\pm},\dots,g_d^{\pm}\}$. Hence,
$$\mu_{\pi,S}(x^{m})=\frac{N(d,m)\chi_{\pi}(id)}{d(\pi)}  +    \frac{1}{d(\pi)}\mathlarger{\sum}_{\text{nontrivial words}}\chi_{\pi}( g_{i_1}\dots g_{i_m}).$$
Since $\chi_{\pi}(id)=d(\pi)$ we obtain
\begin{equation}\label{momentm} \mu_{\pi,S}(x^{m})=N(d,m)  +   \frac{1}{d(\pi)}\mathlarger{\sum}_{\text{nontrivial words}}\chi_{\pi}( g_{i_1}\dots g_{i_m}).\end{equation}
We proceed by computing the $m$th moment of  Kesten-McKay's law with parameter $2d$. Kesten-McKay's law with parameter $2d$ is the spectral measure of the infinite $2d$-regular tree.
 The $m$th moment of the spectral measure is the average number of cycles (we allow backward paths) of length $m$ on a graph.  We note that the $m$th moment of  Kesten-McKay's law $f_{2d}$ is $N(d,m)$. Becuase, there is a one to one correspondence between the cycles of length $m$ starting at the root on the infinite $2d$-regular tree and the set of trivial words with formal letters $\{ g_1^{\pm},\dots,g_d^{\pm}\}$. We proceed by showing that the  probability that a random generator set $S=\{ g_1^{\pm},\dots,g_d^{\pm}\}$ represents $ \pm id$ with a nontrivial word of length $m$ in $S$ where $m \leq (1/2-o(1))\log_{2d-1}(p)$ is $O(p^{-1/4})$.    We cite \cite[Lemma 10]{GHSS} and \cite[Proposition 11]{GHSS}.
 
 \begin{lem}\label{citlem} \normalfont \cite[Lemma 10]{GHSS} 
  Let $\omega$ be a word of length $l$ in the free group $F_k$. If $\omega$ is not identically 1 for every substitution of values from $PGL_2(\mathbb{Z}/p\mathbb{Z})$, then for a random substitution 
  $$Pr[\omega=1]\leq l/p+O(p^{-2}).$$
 \end{lem}

 \begin{lem}\label{citl} \normalfont \cite[Proposition 11]{GHSS}\label{girth}
 For any $k$ the length of the shortest non-trivial word $\omega(x_1,\dots,x_k)$ such that $f_{\omega}(g_1,\dots,g_k)=1$ for all $g_1,\dots,g_k$ in $SL_2(\mathbb{Z}/p\mathbb{Z})$ is at least $\Omega(p/\log(p))$.
 \end{lem}
\noindent Let $M:=\frac{1}{2}\log_{2d-1}(p)$. Since $M<\Omega(p/\log(p))$ then by Lemma~\ref{citl}, any nontrivial word $\omega$ of length less than $M$ in the free group $F_d$ is not identically 1 for every substitution of values from $SL_2(\mathbb{Z}/p\mathbb{Z})$.  By Lemma~\ref{citlem}, the probability of each such word represent identity is at most $M/p+O(p^{-2})$. We note that  the number of nontrivial words with letters in $\{ g_1^{\pm},\dots,g_d^{\pm}\}$ of length less than $M$ is about $(2d-1)^{(M+1)}$.  So by the union bound the probability that the set $S=\{ g_1^{\pm},\dots,g_d^{\pm}\}$ represents $ \pm id$ with a nontrivial word of length bounded by $M$ in $S$ is less than $\frac{M(2d-1)^{M+1}}{p}=O(p^{-1/4})$.  We denote the complement of these generator sets by $\mathcal{F}$. We note that $S\in \mathcal{F}$ with probability $1-O(p^{-1/4})$.  Note that each $g_i$, $1 \leq i \leq d$ has $|G|$ choices. So, there are $|G|^d$ symmetric subsets $S$ and hence $|\mathcal{F}|=(1-O(p^{-1/4}))|G|^d$. 
Next,  we compute the variance of the $m$th moment of  $\mu_{\pi,S}$ as $S:=\{ g_1^{\pm},\dots,g_d^{\pm}\}$  varies in $\mathcal{F}$.  We denote the variance with
 \begin{equation*}
 \text{Var}(m):=\frac{1}{|\mathcal{F}|}\mathlarger{\sum_{S\in \mathcal{F}}} (\mu_{\pi,S}(x^m)-N(d,m))^2.
 \end{equation*}
 We use equation (\ref{momentm}) and obtain
  \begin{equation}
 \text{Var}(m):=\frac{1}{|\mathcal{F}|}\mathlarger{\sum_{S\in \mathcal{F}}} \Big(\frac{1}{d(\pi)}\mathlarger{\sum}_{\text{nontrivial words}}\chi_{\pi}( g_{i_1}\dots g_{i_m})\Big)^2.
  \end{equation}
  We apply Cauchy-Schwarz inequality to obtain
   \begin{equation}\label{varr}
 \text{Var}(m)\leq \frac{1}{|\mathcal{F}| d^2(\pi)} \mathlarger{\sum_{S\in \mathcal{F}}} (2d)^{m}\Big(\mathlarger{\sum}_{\text{nontrivial words}}\chi_{\pi}^2( g_{i_1}\dots g_{i_m})\Big),
  \end{equation}
  where $ (2d)^{m}$ is the number of total words of length $m$  with letters in $S:=\{ g_1^{\pm},\dots,g_d^{\pm}\}$. We give an upper bound on the right hand side of equation (\ref{varr}). 
Let $N(S,m,g)$ denote the number of nontrivial words with length $m$ and letters in $S:=\{ g_1^{\pm},\dots,g_d^{\pm}\}$ that represents $g \in SL_2(\mathbb{Z}/p\mathbb{Z}) $.   
From the definition of $\mathcal{F}$, we have $N(S,m,\pm id)=0$ for every $S\in \mathcal{F}$. Hence,
$$\mathlarger{\sum_{S\in \mathcal{F}}} \mathlarger{\sum}_{\text{nontrivial words}}\chi_{\pi}^2( g_{i_1}\dots g_{i_m})=  \mathlarger{ \sum_{g\neq \pm \rm{id} } }  \sum_{S\in \mathcal{F}}N(S,m,g) \chi^2_{\pi}(g). $$
We give an upper bound on $\sum_{S\in \mathcal{F}}N(S,m,g)$ for every $g\neq\pm\rm{id}$. Note that $SL_2(\mathbb{Z}/p\mathbb{Z}) $ acts by conjugation on $\mathcal{F}$ and that action implies that $N(S,m,g)=N(hSh^{-1},m,hgh^{-1})$ for every  $h\in SL_2(\mathbb{Z}/p\mathbb{Z}) $.  If $g\neq \pm id$ then the conjugacy class  $[g]$ in $SL_2(\mathbb{Z}/p\mathbb{Z})$ has at least  $(p^2-1)/2$ elements; see  \cite[Page 71]{Fulton}. We have $|\mathcal{F}|$ symmetric subsets $S$  and there are less than $(2d)^{m}$ nontrivial words of length $m$ with letters in $S$. We double count the number of the nontrivial words of length $m$ as we vary $S$ and words $g_{i_1}\dots g_{i_m}$:
\begin{equation*}
\begin{split}
 |\mathcal{F}| (2d)^{m} =\mathlarger{\sum_{S\in \mathcal{F}}} (2d)^{m}&\geq  \mathlarger{\sum_{S\in \mathcal{F}}} \sum_{g\in G }  N(S,m,g)
 \\
 &\geq \mathlarger{\sum_{S\in \mathcal{F}}} \sum_{\text {Conj class}[g]\neq [\pm id]  } \frac{p^2-1}{2}  N(S,m,g)
 \\
 &\geq \sum_{\text {Conj class}[g]\neq [\pm id]  } \frac{p^2-1}{2}   \mathlarger{\sum_{S\in \mathcal{F}}}  N(S,m,g).
 \end{split}
\end{equation*}
From the above inequality, we deduce that for every $g\neq \pm id$
\begin{equation}\label{inequal}
\frac{2|\mathcal{F}| (2d)^{m}}{p^2}\geq   \mathlarger{\sum_{S\in \mathcal{F}}} N(S,m,g).
\end{equation}
Finally, we prove an upper bound on $\text{Var}(m)$. Recall from inequality (\ref{varr})
     \begin{equation}
     \begin{split}
 \text{Var}(m)&\leq \frac{1}{|\mathcal{F}| d^2(\pi)} \mathlarger{\sum_{S\in \mathcal{F}}} (2d)^{m}\Big(\mathlarger{\sum}_{\text{nontrivial words}}\chi_{\pi}^2( g_{i_1}\dots g_{i_m})\Big)
 \\
  & =  \frac{(2d)^{m}}{|\mathcal{F}| d^2(\pi)} \mathlarger{ \sum_{g\in G} }  \sum_{S\in \mathcal{F}}N(S,m,g) \chi^2_{\pi}(g).
 \end{split}
  \end{equation}
We use inequality (\ref{inequal}), to obtain 
$$ \text{Var}(m) \leq  \frac{2(2d)^{2m}}{p^2 d^2(\pi)} \sum_{g\in G} \chi^2_{\pi} (g). $$
Since, $\pi$ is a nontrivial  irreducible representation of $SL_2(\mathbb{Z}/p\mathbb{Z})$ then $$\sum_{g\in G} \chi^2_{\pi} (g)=|G|,$$ and the dimension of the irreducible representation $d(\pi)\geq (p-1)/2$. Hence,
 \begin{equation}\label{maineq} \text{Var}(m) \leq \frac{2(2d)^{2m}|G|}{p^2 d^2(\pi)}=O\Big( \frac{(2d)^{2m}}{p}\Big).\end{equation}
Therefore, for  $m \leq (1/4)\log_{2d}(p)$, we obtain
\begin{equation} \label{Vaar}\text{Var}(m)=O(p^{-1/2}).\end{equation}
By Chebyshev's inequality, we deduce that for every $1\leq m \leq (1/4)\log_{2d}(p)$
$$|\mu_{\pi,S}(x^m)-f_{2d}(x^m)|\geq p^{-1/10}. $$
with probability at most $O(p^{-1/4})$. This shows inequality~\ref{mme} holds with probability $1-O(p^{-1/4})$ and completes the proof of Theorem~\ref{ddd}. %
\end{proof}
\noindent
%
Finally, we give a proof of Corollary~{\ref{cccc}}.
 Recall that an eigenvalue $\lambda$ of $\pi (g_1+g^{-1}_1+\dots+ g_{d}+g^{-1}_d)$ is an exceptional eigenvalue if $\lambda \notin [-2\sqrt{2d-1}, 2\sqrt{2d-1}]$ and $N(\alpha, \pi,S)$ is the number of exceptional eigenvalues such that $|\lambda| > \alpha (2\sqrt{2d-1})$ where $\alpha >1$. Note that the number of monochromatic eigenvalues is $d(\pi)= O(p)$. So, the trivial upper bound is $N(\alpha, \pi,S) =O(p)$.  In Corollary~{\ref{cccc}}, we give a sub-exponential  upper bound on $N(\alpha, \pi,S)$.
\begin{proof}[Proof of Corollary~(\ref{cccc})]\normalfont
\noindent  We note that by equation (\ref{Vaar}) 
$$ \label{Var}\text{Var}(m)= \frac{1}{|\mathcal{F}|}\mathlarger{\sum_{S\in \mathcal{F}}} (\mu_{\pi,S}(x^m)-f_d(x^m))^2= o(p^{-\delta}),$$
where $m \leq (1/2-\delta)\log_{2d}(p)$. By Chebyshev's inequality, we deduce  for  every $\varepsilon>0$ and with probability $(1-O(p^{-\delta}))$ for a symmetric subsets  $S$ 
\begin{equation}\label{mmm}
|\mu_{\pi,S}(x^m)- f_d(x^m)|<\varepsilon.
\end{equation} 
We note that the support of Kesten-McKay's  law $f_{2d}$ is inside $[-2\sqrt{2d-1}, 2\sqrt{2d-1}]$. Hence,
\begin{equation}\label{rrr}
f_d(x^m)\leq (2\sqrt{2d-1})^{m}.
\end{equation}  
From inequalities (\ref{mmm}) and (\ref{rrr}), we obtain
$$\mu_{\pi,S}(x^m) \leq  (2\sqrt{2d-1})^{m}.$$
By Markov's inequality and the definition of  $N(\alpha, \pi,S)$, we obtain
$$\frac{N(\alpha, \pi,S)\big(\alpha (2\sqrt{2d-1})\big)^{m}}{p} \leq (2\sqrt{2d-1})^{m}.$$
Hence,
$$N(\alpha, \pi,S)\leq p/\alpha^{m}.$$
We let $m = (1/2-o(1))\log_{2d}(p)$ and deduce that
$$N(\alpha, \pi,S)\leq p^{1-\frac{(\log_{2d}{\alpha})+o(1)}{2} } .$$
This concludes the proof of Corollary~{\ref{cccc}}.
\end{proof}
\begin{rem} \normalfont
If we prove for $m \leq 2\log_{2d}(p) +O(1)$ 
$$\label{Vaar}\text{Var}(m)=o(p^{-\delta}),$$
then $N(\alpha, \sqrt{2d},S)=0 $. Therefore, we have spectral gap for monochromatic eigenvalues. If $\label{Vaar}\text{Var}(m)=o(p^{-\delta})$  for every  $m\gg \log_{2d}(p)$ then monochromatic eigenvalues are almost Ramanujan. 
\end{rem}

\section{The top monochromatic eigenvalue}\label{top}
 In this section, we present our numerical results on the top monochromatic eigenvalue  of Cayley graphs of ${\rm SL}_2[\mathbb{Z}/p\mathbb{Z}]$. Our numerical results for the fixed generator family and random family of Cayley graphs of $SL_2(\mathbb{Z}/p\mathbb{Z})$  show that  the nontrivial eigenvalues are bounded by  Ramanujan bound 
$$|\lambda| < 2(1+\varepsilon)\sqrt{d-1},$$ 
as the prime number $p$ goes to infinity.  Our numerical results together with Corollary~\ref{cccc} suggest that random  and fixed generator family of Cayley graphs of $SL_2(\mathbb{Z}/p\mathbb{Z})$ are almost Ramanujan graphs. This generalizes the celebrated result of Bourgain and Gamburd on the expansion of  random  and fixed generator family of Cayley graphs. Note that our experiments  are consistent with the analogous results for the random regular graphs. By a result of Friedman \cite{free}, random $d$-regular graphs are almost Ramanujan graphs where $d$ is fixed and the number of vertices does to infinity. 

 The spectrum of Cayley graph of $SL_2(\mathbb{Z}/p\mathbb{Z})$ is the union of   principal series and discrete series representations. The top non-trivial eigenvalue of $SL_2(\mathbb{Z}/p\mathbb{Z})$ is the maximum over all monochromatic eigenvalues. As observed by Lafferty and Rockmore \cite{LR3}, different monochromatic eigenvalues  of $SL_2(\mathbb{Z}/p\mathbb{Z})$  have similar statistical  properties.  So, we restrict ourself to the the top eigenvalue of all principal series representations and the  top eigenvalue of Steinberg representation.   We present  our numerical results for the top nontrivial eigenvalue of affine Cayley graphs that is associated to the top eigenvalue of all principal series representations. We also give our numerical results for the top nontrivial eigenvalue of projective Cayley graph that is associated to Steinberg representation. Based on our experiments, the top monochromatic eigenvalue is asymptotic to  Ramanujan bound  $2\sqrt{d-1}$ as $p$ goes to infinity. 
 As pointed out in the introduction, different monochromatic  eigenvalues  are independent of each other, since the number of principal series representations  is linear in $p$  and each oscillates around the Ramanujan bound, then heuristically the maximum of them would be bigger than the Ramanuajn bound. We verify numerically that the top nontrivial eigenvalue of  affine Cayley graphs is always bigger than the Ramanujan bound while the top eigenvalue of the projective Cayley graphs is less than the Ramanujan bound with a positive probability as prime number $p$ goes to infinity. 
\subsection{Fixed generator families of Cayley graphs of ${\rm SL}_2[\mathbb{Z}/p\mathbb{Z}]$}

In this section, we present our numerical experiments for two families of Schreier graphs obtained by  the fixed generator set $$S:=\Big\{\begin{bmatrix} 1&\pm2 \\ 0 &1     \end{bmatrix}, \begin{bmatrix} 1&0 \\ \pm2 &1     \end{bmatrix}   \Big\}.$$ The two families are projective Cayley graph of $\mathbb{P}^1(\mathbb{Z}/p\mathbb{Z})$, and affine Cayley graphs of $\mathbb{A}^2-\{0\}(\mathbb{Z}/p\mathbb{Z})$ (the affine space minus the origin)  generated by the action of  $S$. We  give our numerical experiments on  the top non-trivial eigenvalue of these graphs.      We numerically show that the top eigenvalue is asymptotically optimal as $p$ goes to infinity.
 For the Projective Cayley graphs, we  pick a large prime of size $p\sim2000000$, for the affine Cayley graphs $p\sim 3500$. Since,  the order of vertices is $p+1$ comparing to $p^2-1$ vertices for the affine Cayley graphs.

\begin{figure}[t]
\vspace{0.45cm}
\centering
\raisebox{-0.6cm}{
\includegraphics[width=.54\textwidth]{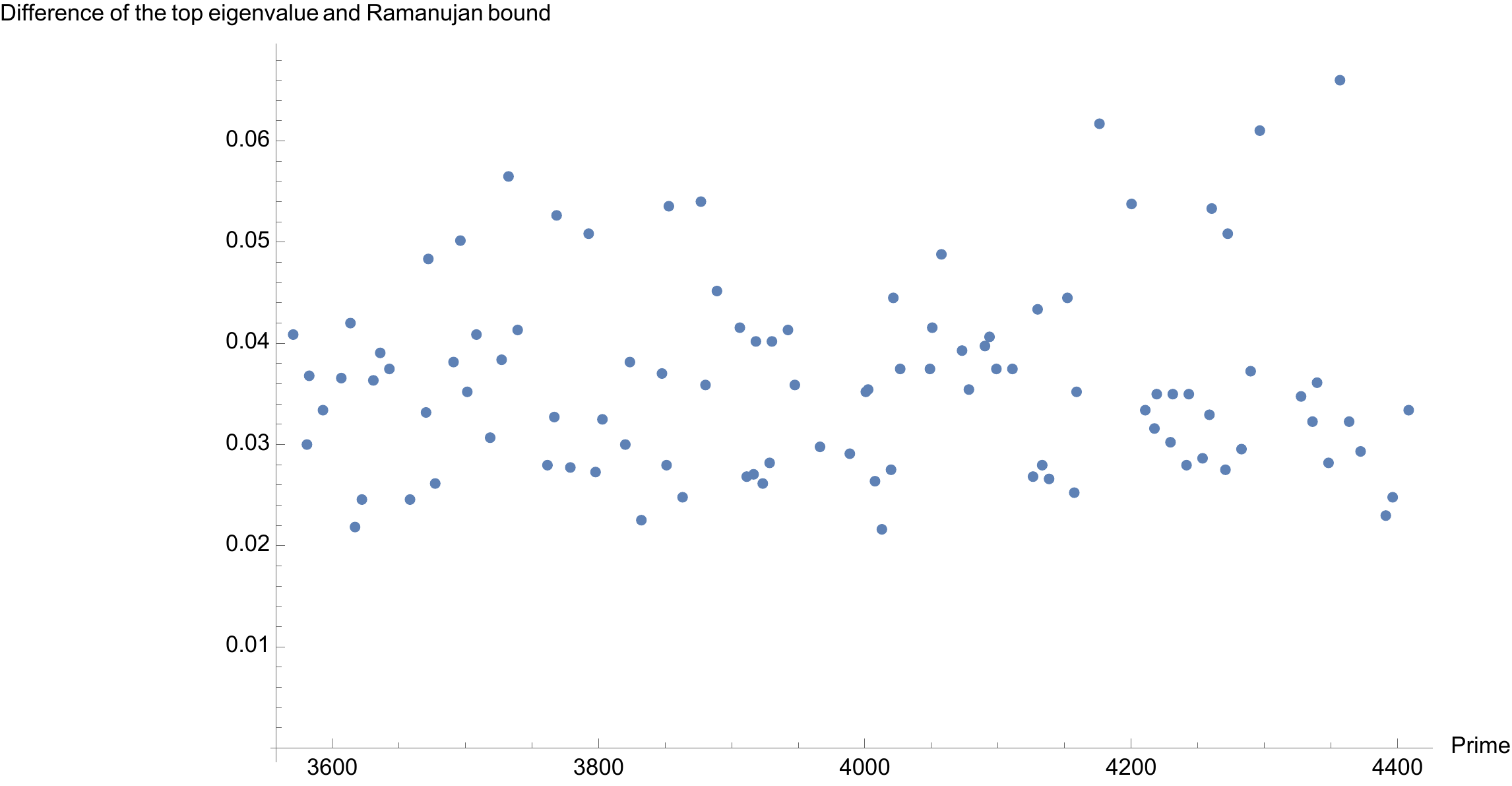}
\includegraphics[width=.54\textwidth]{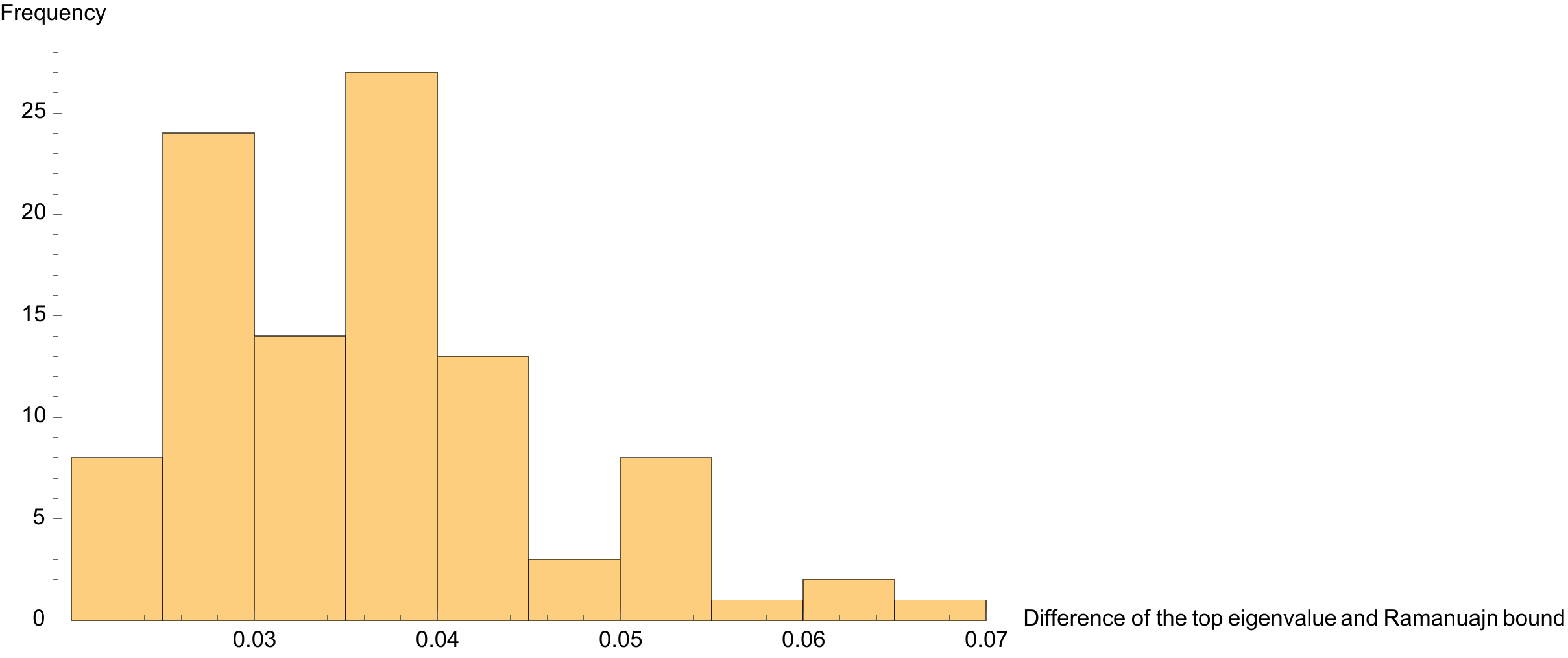}
}\hspace{0.25cm}
\vspace{-0.48cm}
\caption{Difference between the top eigenvalue and Ramanujan bound of affine Cayley graph of $\mathbb{A}^2-\{0\}(\mathbb{Z}/p\mathbb{Z})$ generated by $S$.}
\label{flg}
\vspace{0.2cm}
\end{figure}

 We begin by presenting our experiments for the 4-regular affine Cayley graphs generated by $S$. In Figure~(\ref{flg}), we give the difference between the top non-trivial eigenvalue of these graphs and Ramanujan bound $2\sqrt{3}$ for 100 primes between the 500th prime 3571 and the 600th prime 4409.  The mean and  variance is $(0.0361, 0.0094)$, respectively. None of the affine Cayley graphs in this range  is a Ramanujan. However,  they are almost Ramanujan; i.e. the difference between the top non-trivial eigenvalue and Ramanujan bound is small.

\begin{figure}[t]
\vspace{0.45cm}
\centering
\raisebox{-0.6cm}{
\includegraphics[width=.54\textwidth]{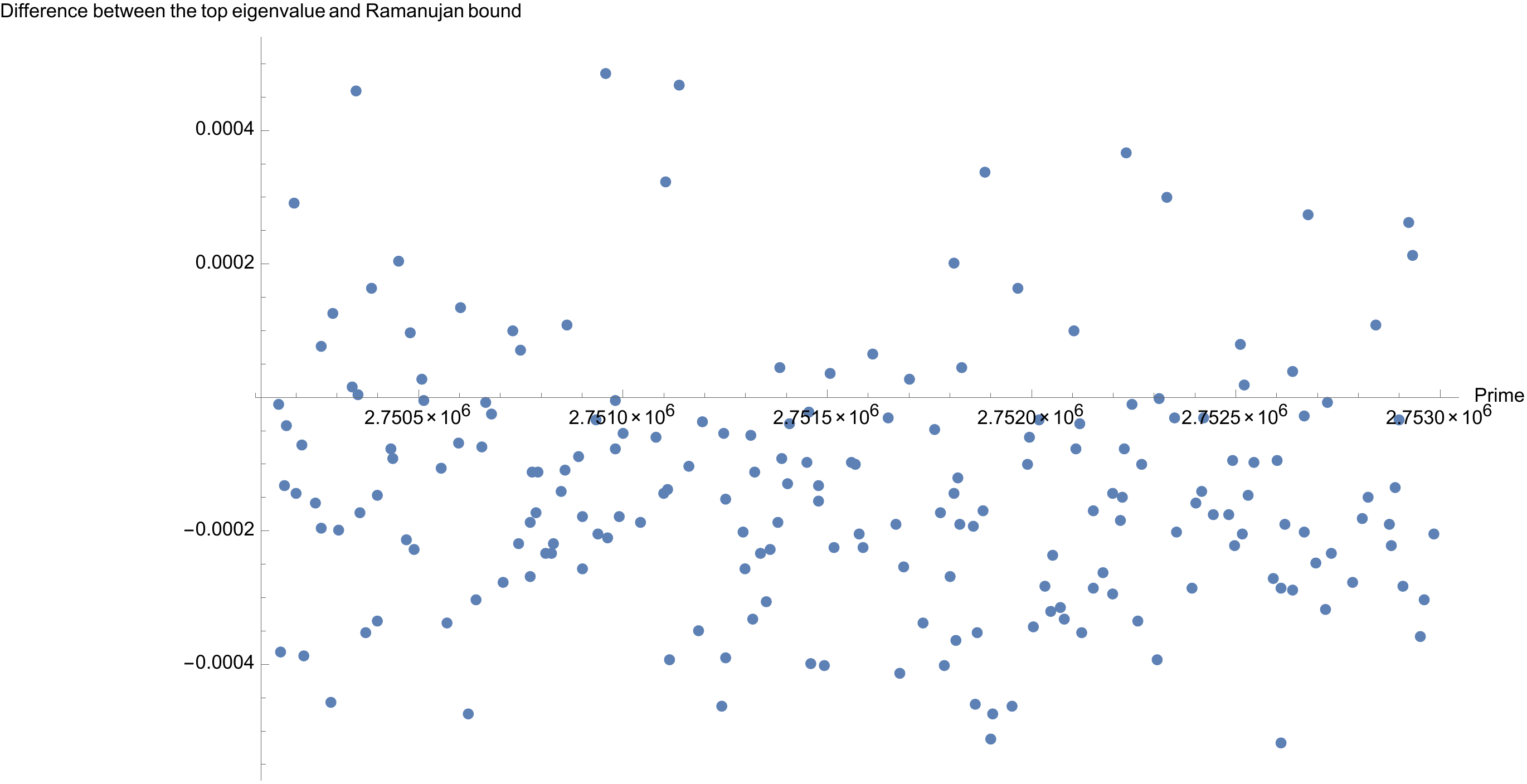}
\includegraphics[width=.54\textwidth]{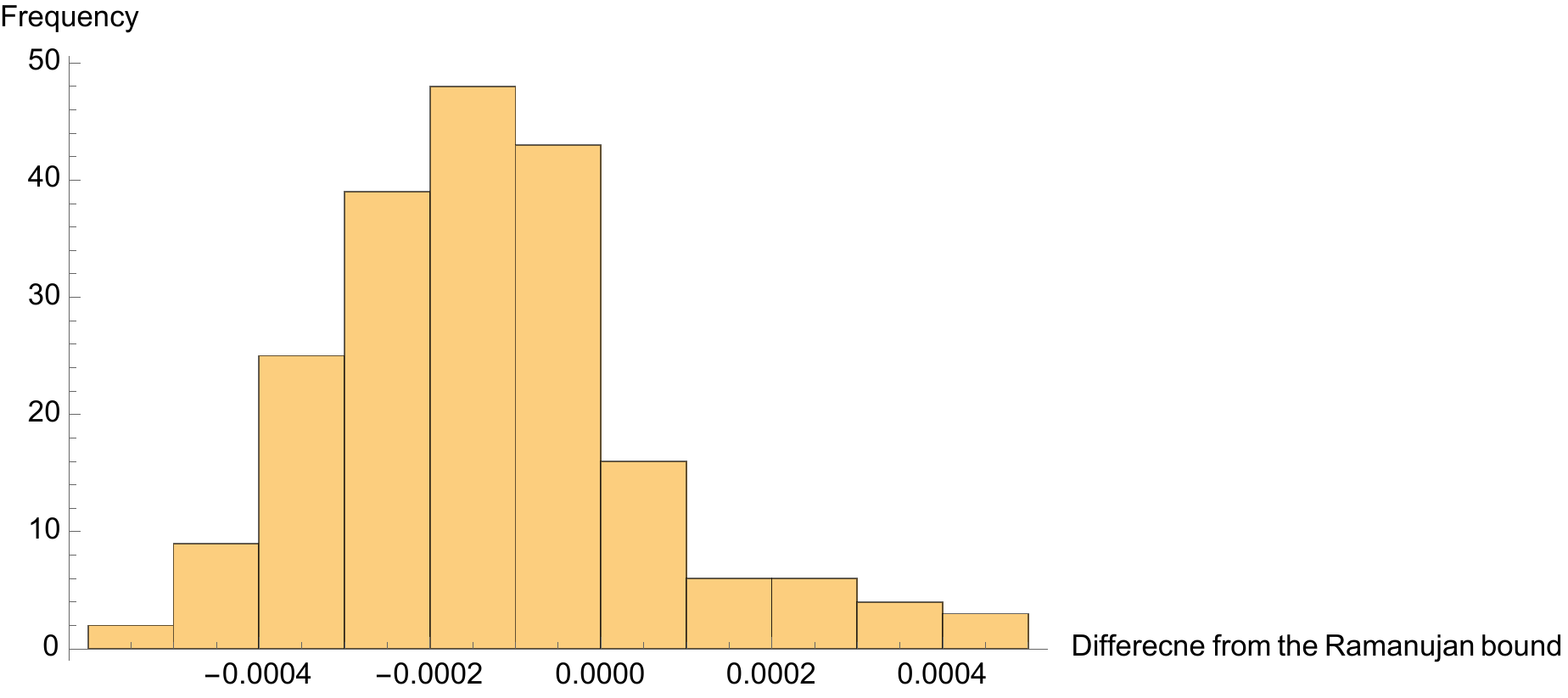}
}\hspace{0.25cm}
\vspace{-0.48cm}
\caption{Difference between the top eigenvalue and Ramanujan bound of the projective Cayley graph of $\mathbb{P}^1(\mathbb{Z}/p\mathbb{Z})$ generated by $S$.}
\label{fp}
\vspace{0.2cm}
\end{figure}

 Next, we give our numerical results for the 4-regular  projective Cayley graphs generated by $S$. In Figure~(\ref{fp}), we give the difference between the top non-trivial eigenvalue of these graphs and Ramanujan bound $2\sqrt{3}$ for 200 primes between the 200000th prime 2750159 and the 200200th prime 2752987.  The mean and the variance is $(-0.000134809, 0.0001875)$, respectively. Note that the mean is negative and only 35 out of 200 projective Cayley graphs in this range are non-Ramanujan.

\subsection{Random families of Cayley graphs}

In this section, we present similar numerical experiments for two random families of  Schreier graphs of ${\rm SL}_2[\mathbb{Z}/p\mathbb{Z}]$.   The two families are  projective Cayley graph of $\mathbb{P}^1(\mathbb{Z}/p\mathbb{Z})$, and affine Cayley graph of $\mathbb{A}^2-\{0\}(\mathbb{Z}/p\mathbb{Z})$ generated by the action of   $R=\{g_1,g_1^{-1},g_2,g_2^{-1} \}$ where $\{g_1,g_2\}$ are chosen uniformly among the pairs of generators. For the affine Cayley graphs, we fix prime number $p= 3571$ and pick $100$ samples of random generator sets $R$ of ${\rm SL}_2[\mathbb{Z}/3571\mathbb{Z}]$. Similarly, for the projective Cayley graphs, we  fix  prime number $p=2750159$ and pick  $200$ samples of random generator sets $R$ of ${\rm SL}_2[\mathbb{Z}/2750159\mathbb{Z}]$. 
 
\begin{figure}[t]
\vspace{0.45cm}
\centering
\raisebox{-0.6cm}{
\includegraphics[width=.54\textwidth]{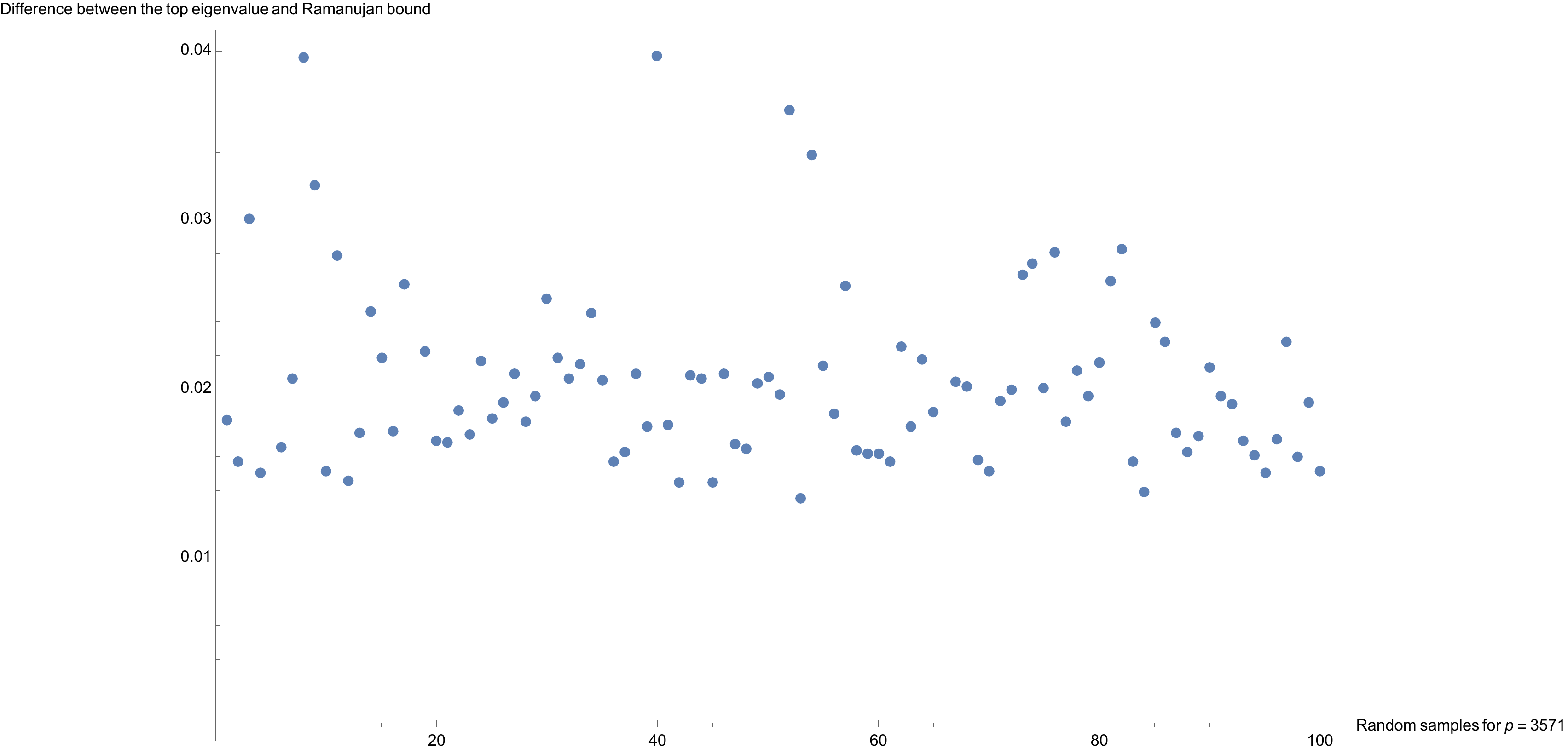}
\includegraphics[width=.54\textwidth]{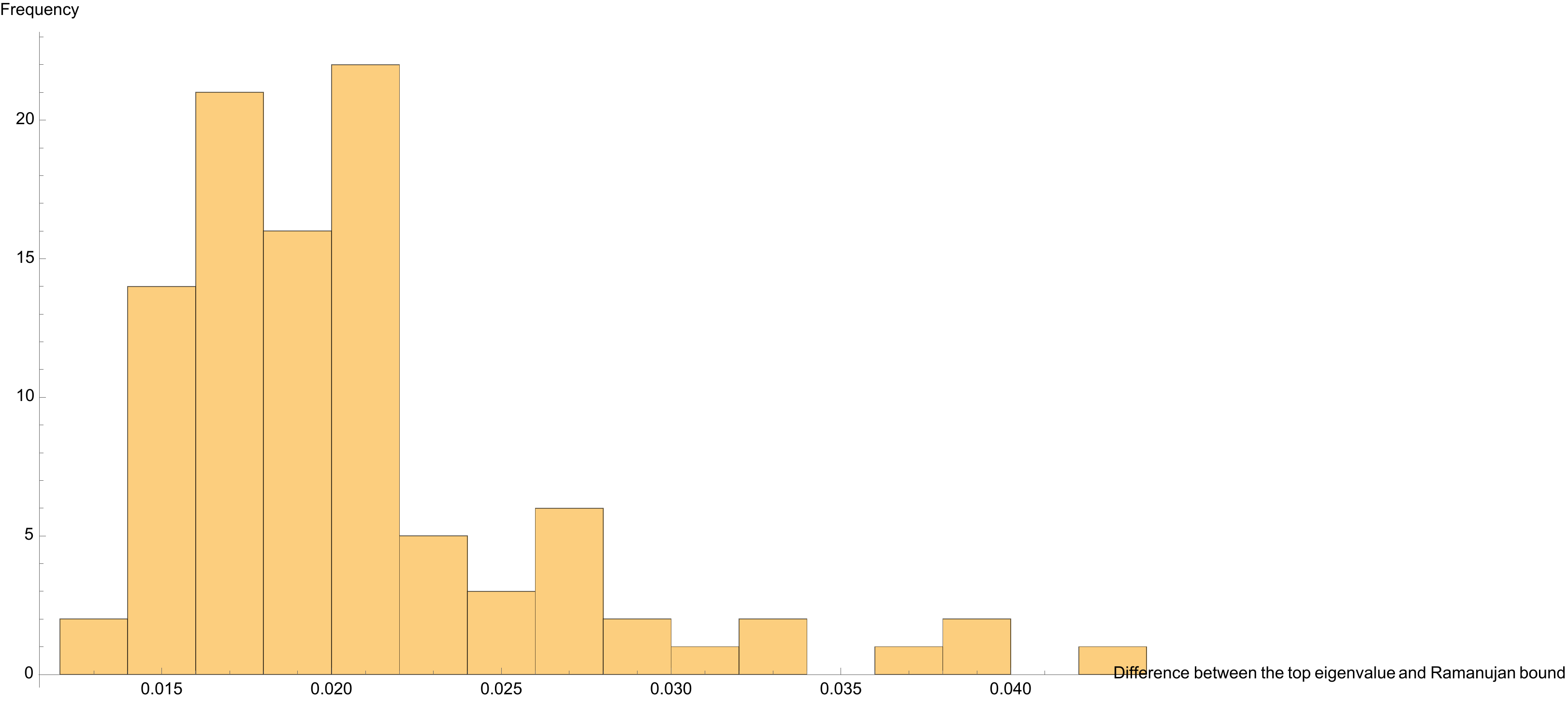}
}\hspace{0.25cm}
\vspace{-0.48cm}
\caption{Difference between the top eigenvalue and Ramanujan bound of 100 samples of the random affine Cayley  4-regular graph for $p=3571$  .}
\label{rl}
\vspace{0.2cm}
\end{figure}

 In Figure~(\ref{rl}), we give the difference between the top non-trivial eigenvalue of random affine graphs and Ramanujan bound $2\sqrt{3}$ for 100 random samples of generator sets of  ${\rm SL}_2[\mathbb{Z}/3571\mathbb{Z}]$.  The mean and  variance is $(0.0225769, 0.0148922)$, respectively. Similar to the affine Cayley graphs generated by the fixed generator  set $S$, none of them  is a Ramanujan. However,  they are almost Ramanujan; i.e. the difference between the top non-trivial eigenvalue and Ramanujan bound is small.  This is consistent with our conjecture about the independence of eigenvalues for different irreducible representation. 
 
\begin{figure}[t]
\vspace{0.45cm}
\centering
\raisebox{-0.6cm}{
\includegraphics[width=.54\textwidth]{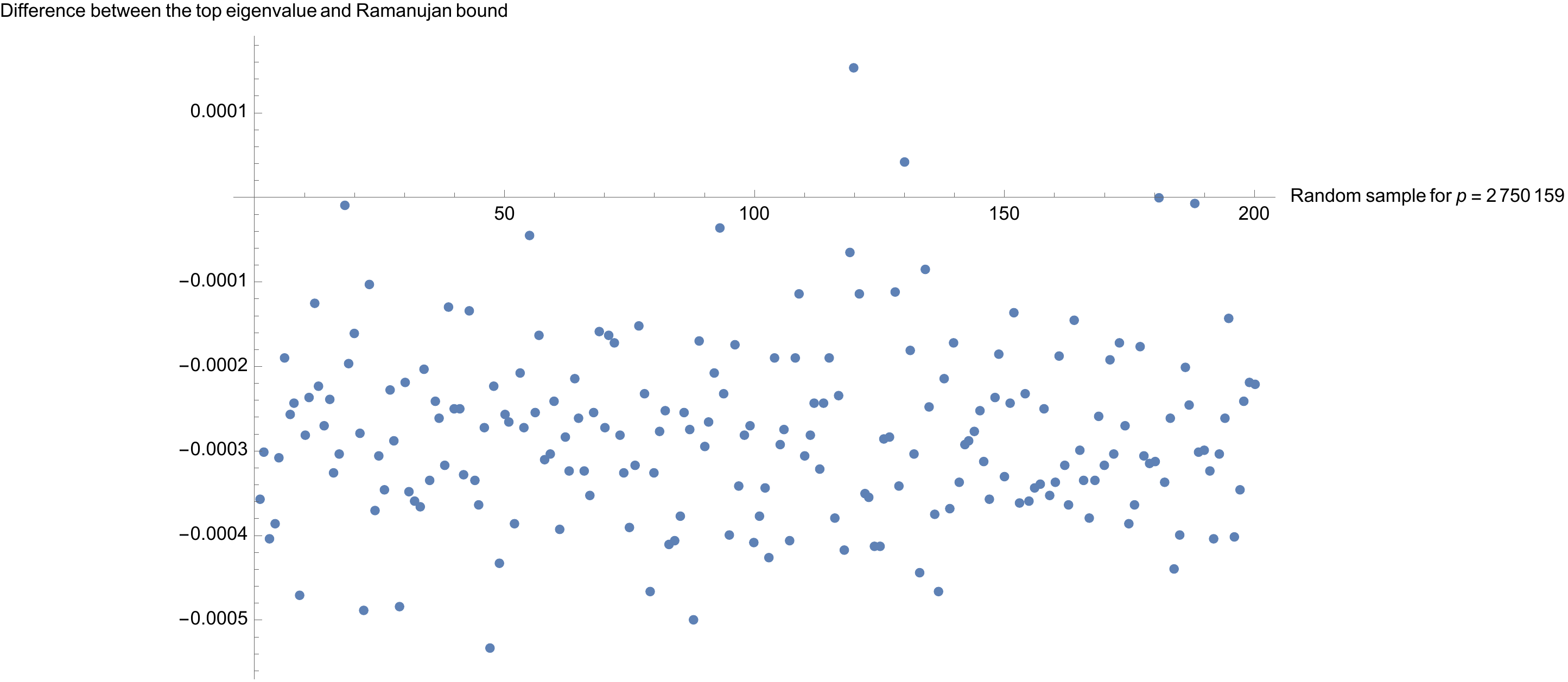}
\includegraphics[width=.54\textwidth]{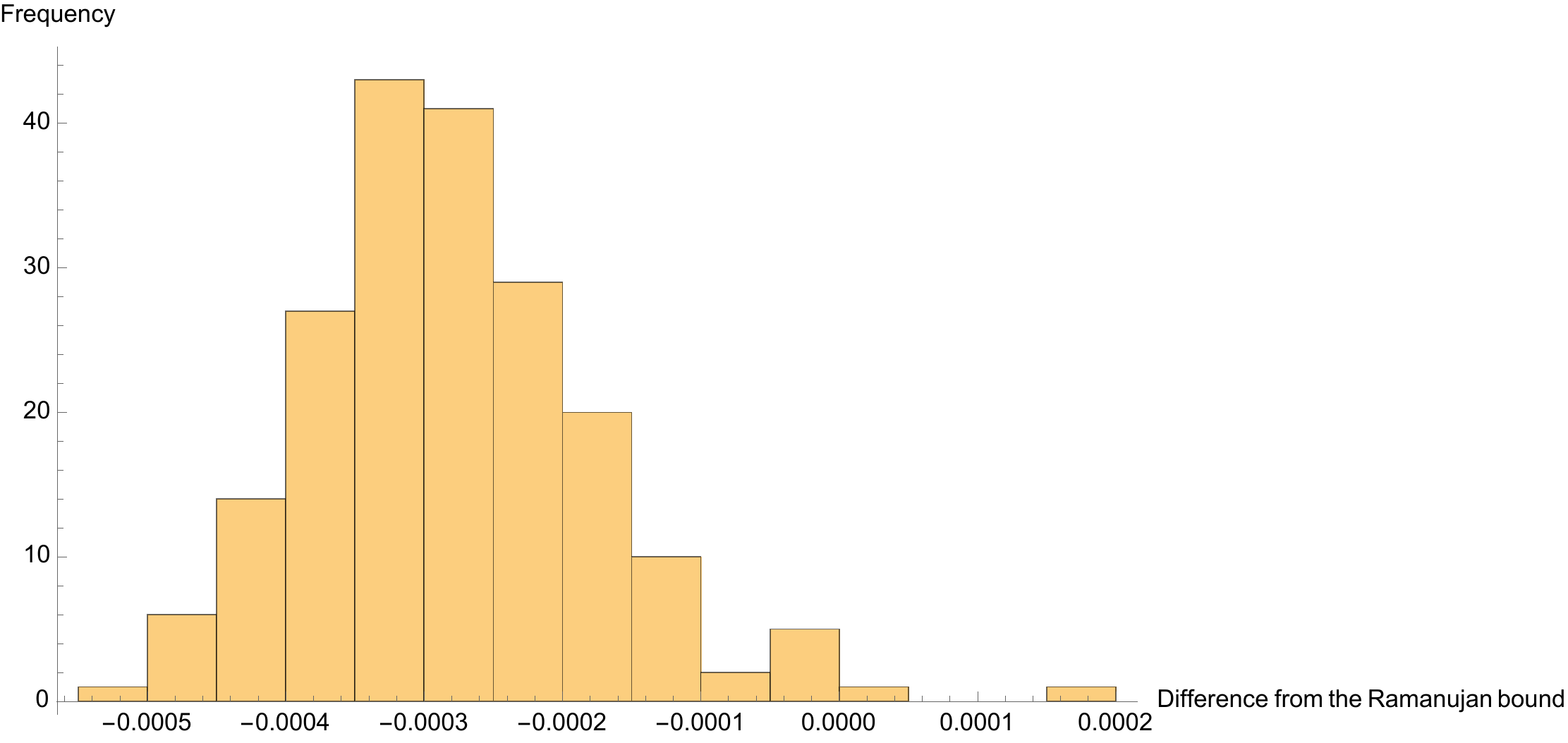}
}\hspace{0.25cm}
\vspace{-0.48cm}
\caption{Difference between the top eigenvalue and Ramanujan bound of 200 samples of the random projective Cayley 4-regular  graph for $p=2750159$.}
\label{rp}
\vspace{0.2cm}
\end{figure}

Finally in Figure~(\ref{rp}),  we give the difference between the top non-trivial eigenvalue of random projective graphs and Ramanujan bound $2\sqrt{3}$ for 200 random samples of generator sets of  ${\rm SL}_2[\mathbb{Z}/2750159\mathbb{Z}]$.  The mean and the variance is $(-0.000279242, 0.000103791)$, respectively. Note that the mean is negative and only 2 out of 200 random samples of projective Cayley graphs in this range are non-Ramanujan.

\section{Diameter}\label{diam}

 In this section, we present our numerical experiments on the diameter of the fixed and random generator families of Cayley and projective Cayley  graphs. Our experiments  show that random generator families of Cayley and projective Cayley graphs have optimal diameter. On the other hand, the diameter of the Cayley graph with LPS generators that is a special fixed generator family of Cayley graphs is not optimal. The second name author in  \cite{Naser}, showed that the diameter of a  family  of LPS Ramanujan graphs is greater than $4/3\log_{d-1}(n)$  where $n$ is the number of vertices and $d$ is the degree of the graph. This lower bound is related to the repulsion of integral points lying on quadrics. The second name author formulated a conjecture on the optimal exponent for the approximation of points on 3 dimensional quadrics   by integral points lying on the quadric  \cite[Conjecture 1.9]{Optimal}. That conjecture  implies that the diameter of LPS Ramanujan graphs is asymptotic to $4/3\log_{d-1}(n)$. More recently,  he developed and implemented a polynomial time algorithm for navigating LPS Ramanujan graphs under a polynomial time algorithm for factoring integers and an arithmetic conjecture on the distribution of numbers representable as a sum of two squares \cite{Sss}. The numerical results from that algorithm strongly support that the diameter of LPS Ramanujan graphs is asymptotic to $4/3\log_{d-1}(n)$; see \cite{Sss}.   This shows that the diameter of Random Cayley graphs  is asymptotically smaller than the diameter of  LPS Cayley graphs. In this paper, we present new numerical results for the diameter of the LPS projective Cayley graphs.  Our numerical results show that the diameter of LPS projective Cayley graphs is optimal.  

\subsection{Diameter of Cayley graphs of ${\rm SL}_2(\mathbb{Z}/p\mathbb{Z})$}
\begin{figure}[t]
\vspace{0.45cm}
\centering
\raisebox{0.2cm}{
\includegraphics[width=.54\textwidth]{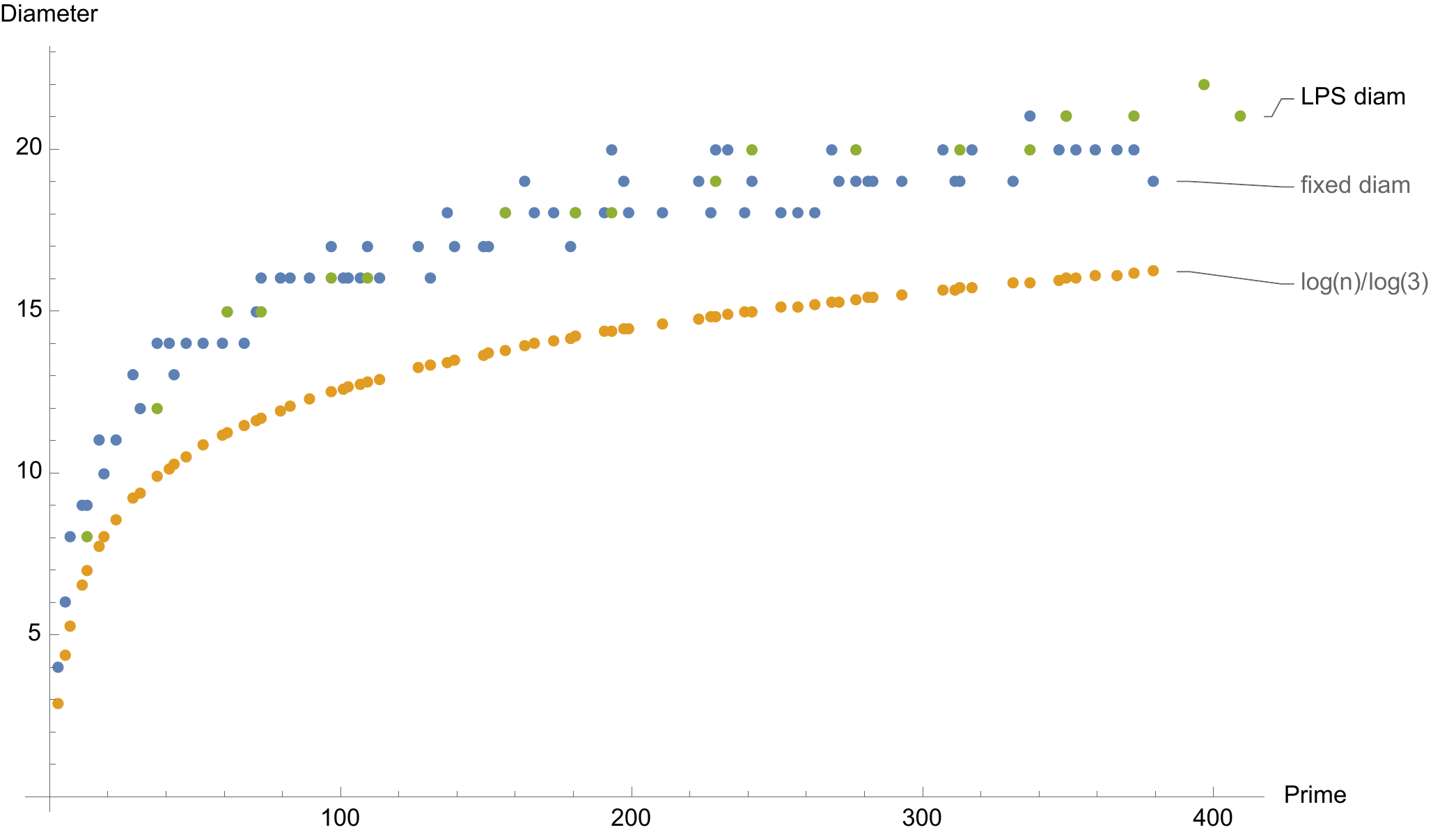}
\includegraphics[width=.54\textwidth]{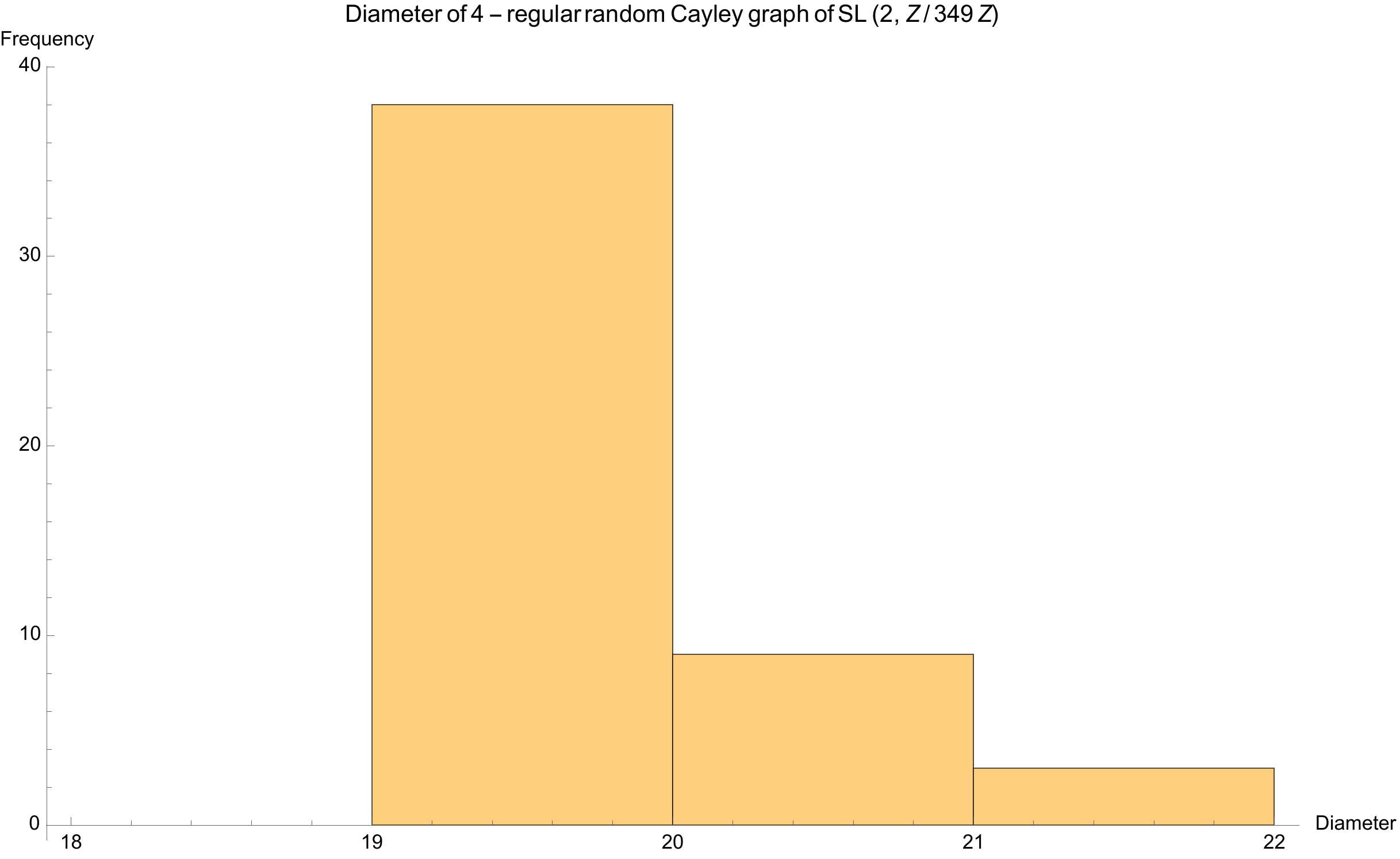}
}\hspace{0.25cm}
\vspace{-0.48cm}
\caption{Diameter of 4-regular Cayley graph of $SL_2(\mathbb{Z}/p\mathbb{Z})$.}
\label{fcresults}
\vspace{0.2cm}
\end{figure}

In this section, we give our numerical results on the diameter of  Cayley graph of  ${\rm SL}_2(\mathbb{Z}/p\mathbb{Z})$ with various symmetric generator sets with 4 elements. We consider Cayley graphs of  ${\rm SL}_2(\mathbb{Z}/p\mathbb{Z})$ generated by the fixed set $S$ where
$$S:=\Big\{\begin{bmatrix} 1&\pm2 \\ 0 &1     \end{bmatrix}, \begin{bmatrix} 1&0 \\ \pm2 &1     \end{bmatrix}   \Big\},$$
and $p$ is a large prime number. Note that the subgroup generated by $S$ is a finite index subgroup of  ${\rm SL}_2(\mathbb{Z})$. Next,   we describe  the LPS generator set associated to prime 3 that has 4 elements. We assume that 3 and $-1$ are quadratic residues mod $p$. This is equivalent to $p \equiv 1 \text{ mod } 12$. We denote the LPS generator set associated to prime 3 by 
$$L:= 1/\sqrt{3}\Big\{ \begin{bmatrix} i & 1 + i \\  -1+i & -i  \end{bmatrix}, \begin{bmatrix}i & -1+i \\ 1+i & -i    \end{bmatrix}, \begin{bmatrix} i & 1-i \\ -1-i & -i   \end{bmatrix},\begin{bmatrix} i & -1-i \\ 1-i & -i  \end{bmatrix} \Big\},$$
where $i=\sqrt{-1}$ mod $p$. The Cayley graph of $SL_2(\mathbb{Z}/p\mathbb{Z})$ generated by $L$ is a Ramanujan graph \cite{Chiu1992}. We also consider the diameter of Cayley graphs of ${\rm SL}_2(\mathbb{Z}/349\mathbb{Z})$  generated by a random symmetric generator set of size 4. Our numerical results for the diameter of these 4-regular Cayley graph of ${\rm SL}_2(\mathbb{Z}/p\mathbb{Z})$  is presented in Figure~\ref{fcresults}. The value of diameter is compered with the trivial lower bound $\log_{3}(n)$, where $n=q^3-q$ is the number of vertices of  Cayley graph of ${\rm SL}_2(\mathbb{Z}/q\mathbb{Z})$. Note that the diameter of the Cayley graph of ${\rm SL}_2(\mathbb{Z}/379\mathbb{Z})$ generated by $S$ is 19. By comparing this to the trivial lower bound $\log_3(n)=16.2137$ and the lower bound for  LPS Ramanujan graphs $4/3\log_3(n)=21.6183$,  we expect that the diameter of the family of Cayley graphs generated by $S$ is asymptotically smaller that the 4 -regular LPS Ramanujan graphs and even optimal.  In Figure~(\ref{fcresults}), we also give the diameter of 50 random samples of Cayley graphs of ${\rm SL}_2(\mathbb{Z}/349\mathbb{Z})$ in the histogram diagram. In our experiments,  38 times the diameter is 19, 9 times is 20 and 3 times is 21. This  shows that the diameter of random Cayley graph of ${\rm SL}_2(\mathbb{Z}/349\mathbb{Z})$ is shorter than the LPS Ramanujan graph. For fixed $d$, we expect that the diameter of $2d$-regular Cayley graphs of ${\rm SL}_2(\mathbb{Z}/p\mathbb{Z})$ generated by a symmetric random set of size $2d$ is optimal and is asymptotic to $\log_{2d-1}(n)$ almost surely as $p$ goes to infinity.

\subsection{Diameter of projective Cayley graph}

\begin{figure}[t]
\vspace{0.45cm}
\centering
\raisebox{0.2cm}{
\includegraphics[width=.54\textwidth]{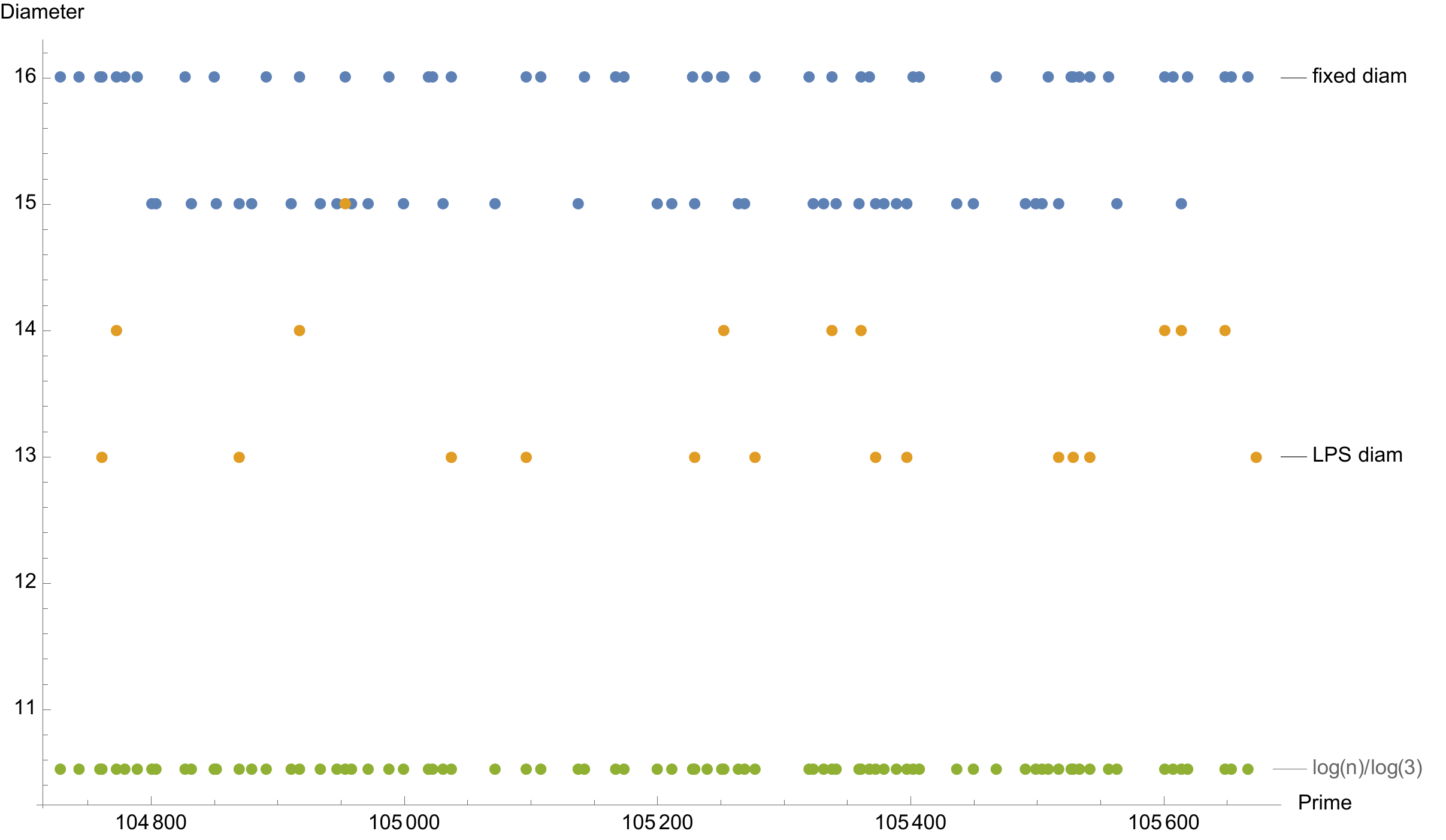}
\includegraphics[width=.54\textwidth]{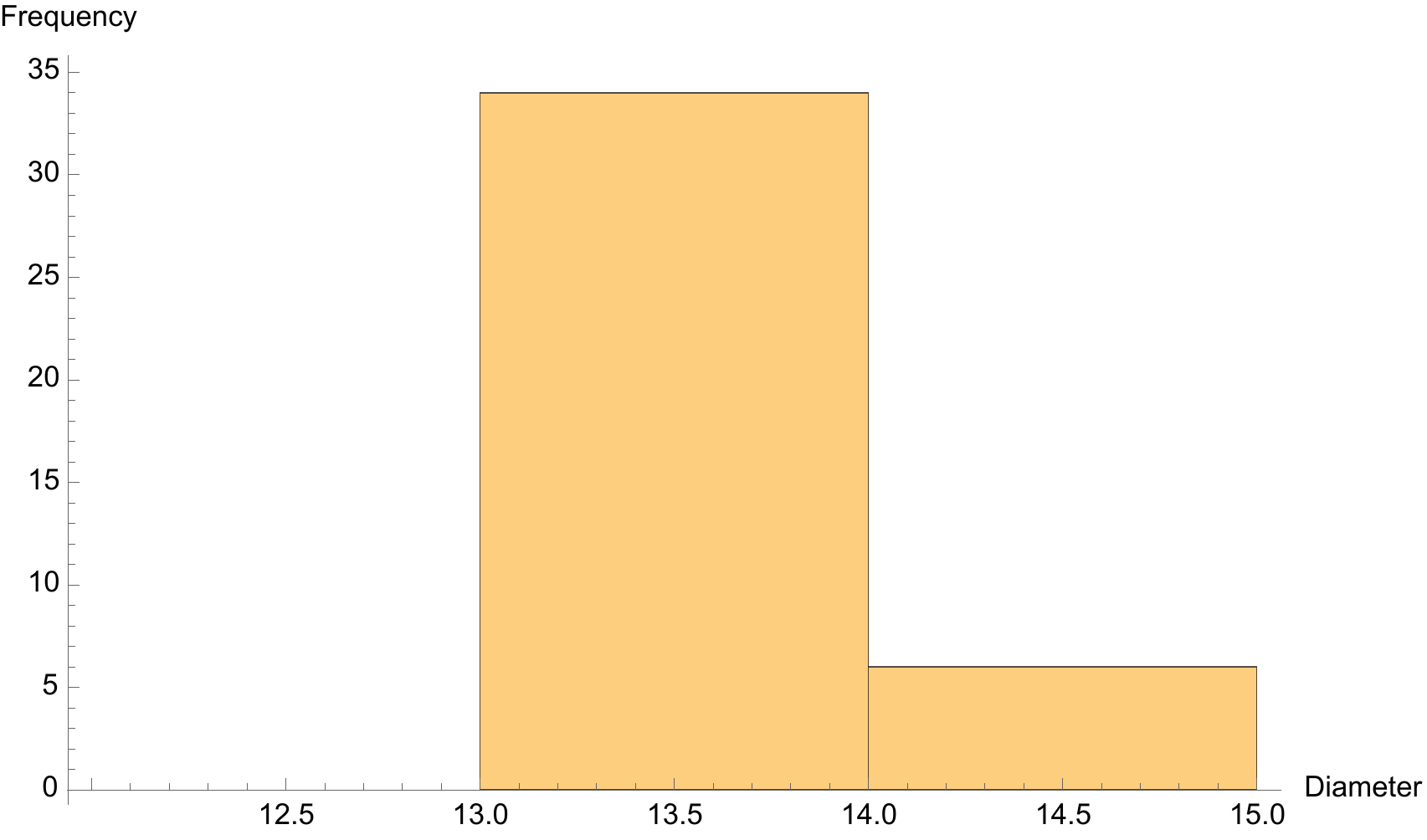}
}\hspace{0.25cm}
\vspace{-0.48cm}
\caption{Diameter of 4-regular projective Cayley graph of $\mathbb{P}^1(\mathbb{Z}/p\mathbb{Z})$.}
\label{diamproj}
\vspace{0.2cm}
\end{figure}

In Figure~(\ref{diamproj}), we give our numerical experiments for the diameter of the projective Cayley graphs of $\mathbb{P}^1({\mathbb{Z}/p\mathbb{Z}})$ of various  generator sets with 4 elements. Since, the Cayley graph of $SL_2(\mathbb{Z}/p\mathbb{Z})$ generated by $L$ is a Ramanujan graph \cite{Chiu1992}. As a result  projective Cayley graph of $\mathbb{P}^1({\mathbb{Z}/p\mathbb{Z}})$ generated by $L$ is also a Ramanujan graph. Unlike the case of LPS Cayley graphs of $SL_2(\mathbb{Z}/p\mathbb{Z})$, the diameter of LPS projective Cayley graphs of $\mathbb{P}^1({\mathbb{Z}/p\mathbb{Z}})$ is optimal. In Figure~(\ref{diamproj}), we compare the diameter of Projective Cayley graph generated by $S$ that is labeled by "fixed diam" with the diameter of LPS projective graphs that is labeled by "LPS diam".  We also give the diameter of 40 random samples of  4-regular projective Cayley graph of $\mathbb{P}^1({\mathbb{Z}/104729\mathbb{Z}})$. 34 times the diameter of random projective graph defined on $\mathbb{P}^1({\mathbb{Z}/104729\mathbb{Z}})$ is 13 and 6 times is 14. In contrast to the previous section, the diameter of LPS projective Cayley graphs is shorter than the diameter of the projective graphs generated by $S$. For prime $p=105541$, the diameter of the projective Cayley graph defined on $\mathbb{P}^1({\mathbb{Z}/105541\mathbb{Z}})$ with LPS generators is $13$ whereas with fixed generator set $S$ is 16. This should be compared with the trivial lower bound $\log_{d-1}n=10.53$ .



\subsection{Radius of the projective Cayley graph centered at zero}
\begin{figure}[t]
\vspace{0.45cm}
\centering
\raisebox{0.2cm}{
\includegraphics[width=.54\textwidth]{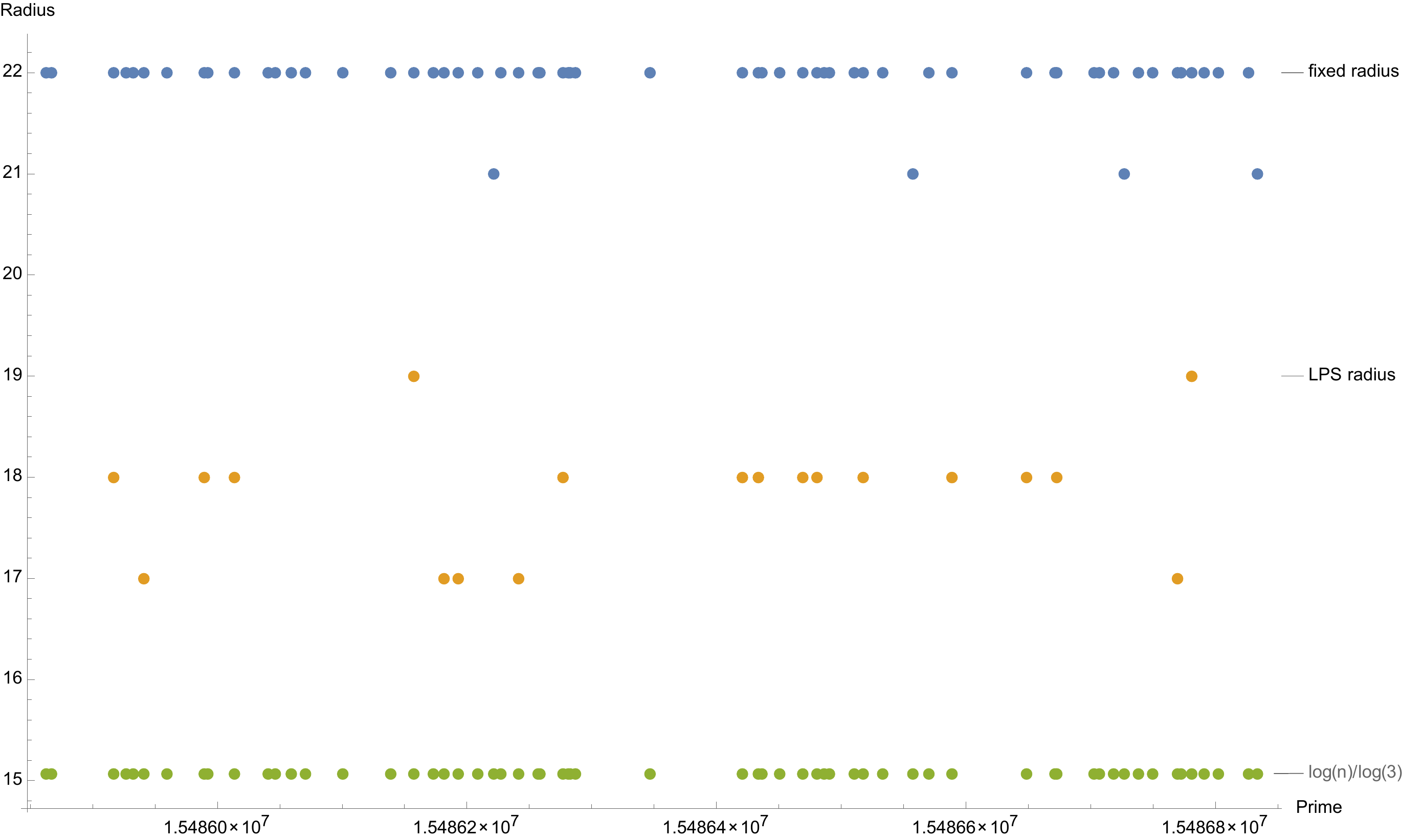}
\includegraphics[width=.54\textwidth]{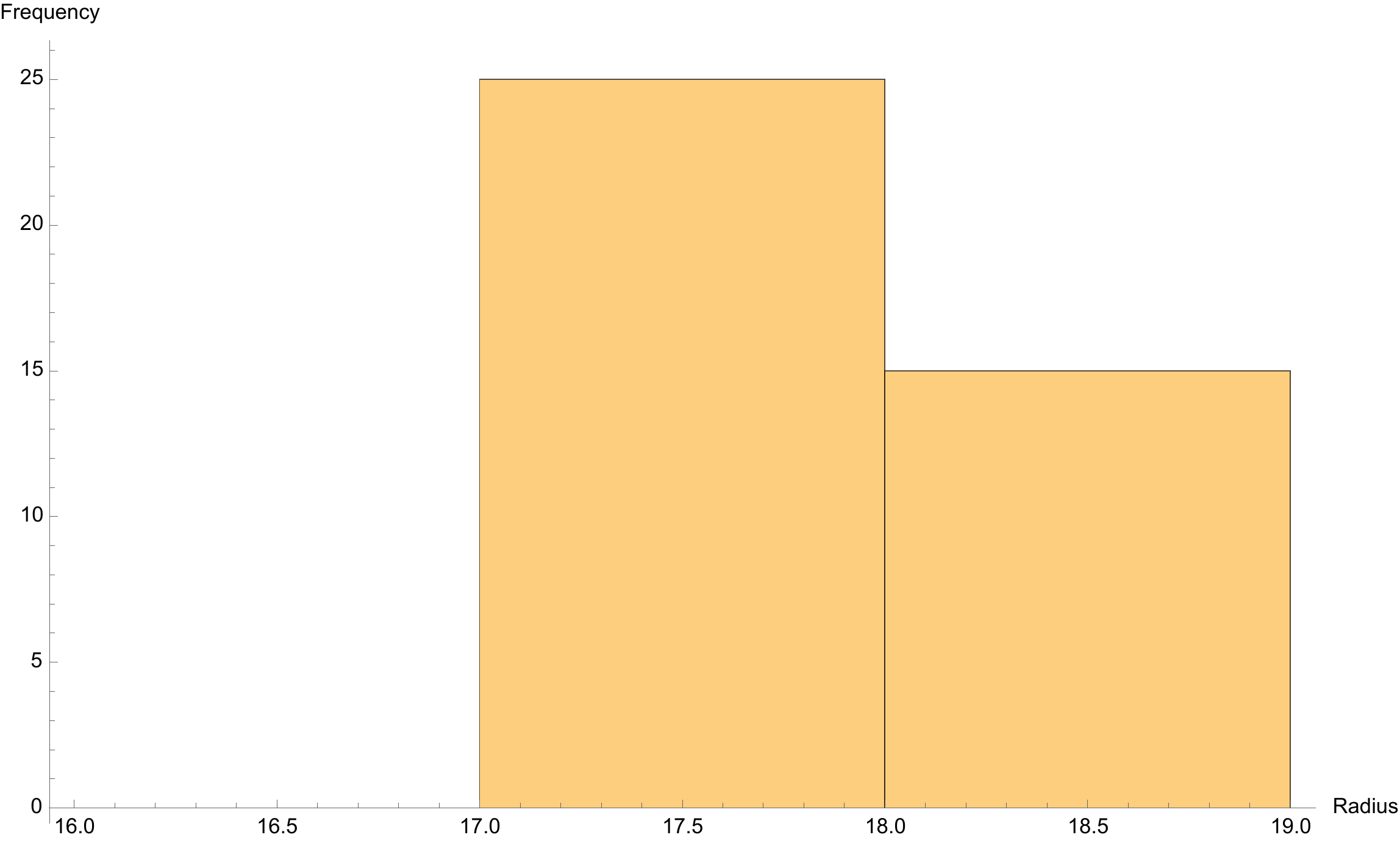}
}\hspace{0.25cm}
\vspace{-0.48cm}
\caption{Radius of the projective graphs of $\mathbb{P}^1(\mathbb{Z}/p\mathbb{Z})$ centered at zero.}
\label{radproj}
\vspace{0.2cm}
\end{figure}

Finally, we give our numerical results for the radius of the projective Cayley graph centered at zero. The size of prime numbers that we experiment with is about $1.5\times 10^7.$ Our numerical results for the radius show that the radius of LPS projective Cayley graphs of $\mathbb{P}^1({\mathbb{Z}/p\mathbb{Z}})$ is optimal.  For prime $p=15486769$ the radius  of LPS projective Cayley graphs is $17$ compered with the trivial lower bound $\log_3(p)=15.1$.  In Figure~(\ref{radproj}), we compare the radius of Projective Cayley graph generated by fixed generator set $S$ and LPS generator set $L$ with the the trivial lower bound $\log_3(n)$ on the radius. It is evident from Figure~(\ref{radproj}), that the radius of LPS projective graphs is smaller than the radius of the fixed generator family.  We also give the radius  of 40 random samples of  4-regular projective Cayley graph of $\mathbb{P}^1({\mathbb{Z}/p\mathbb{Z}})$ for $p=15485863$. 25 times the radius of random projective graph on $\mathbb{P}^1({\mathbb{Z}/p\mathbb{Z}})$ is 17 and 15 times is 18.  

%
%
%
%
%
%

%

\section{Extending the results to $SL_n(\mathbb{Z}/p\mathbb{Z})$}\label{extend}
In this section we prove that the distribution of the monochromatic eigenvalues of $d$-regular random Cayley graphs of $G=SL_n(\mathbb{Z}/p\mathbb{Z})$ for any $n\geq 3$  converges to  Kesten-Mackay law for almost all random generators as $p$ goes to infinity.

\begin{proof}[Proof of Theorem~\ref{lasttheorem}]
We follow the same steps as in the proof of Theorem~\ref{ddd}  for  $G=SL_2(\mathbb{Z}/p\mathbb{Z})$. We reduce bounding the discrepancy  $D(\mu_{\pi,S}, f_{2d})$ to bounding  $|\mu_{\pi,S}(x^m)- f_d(x^m)|$ where $\mu_{\pi,S}(x^m)$ and $f_{2d}(x^m)$  are the  the $m$th moment of $\mu_{\pi,S}$ and $f_{2d}$, respectively. In particular, if $S$ is a generator set for $SL_n[\mathbb{Z}/p\mathbb{Z}]$ such that for every integer  $1 \leq m \leq A$, 
 \begin{equation}\label{mme1}
|\mu_{\pi,S}(x^m)- f_d(x^m)|\ll p^{-\delta},
\end{equation} 
 for some  $0< \delta$.  Then 
 \begin{equation}\label{discp1}D(\mu_{\pi,S}, f_{2d}) \ll 1/A.\end{equation}
Therefore, it is enough to prove inequality $(\ref{mme1})$ holds for generator set $S$ with probability $1-O(p^{-1})$. Our strategy to prove inequality (\ref{mme1}) is to exclude a measure $O(p^{-1})$ family of symmetric subsets $S$ and then give an upper bound  on the variance of the $m$th moment of the remaining symmetric sets $S$. By Chebyshev's inequality, we deduce inequality (\ref{mme1}) for symmetric subsets $S$  with probability $(1-O(p^{-1}))$.  
  In the case of $SL_2(\mathbb{Z}/p\mathbb{Z})$, we excluded the set of all generators (bad generators) that can spell identity with a word of small length. Whereas for $G=SL_n(\mathbb{Z}/p\mathbb{Z})$, we exclude the set  of symmetric generators of size $2d$ that can spell a non-regular element  $g\in G$ with a word of length smaller than $A$.  The regular elements of $G=SL_n(\mathbb{Z}/p\mathbb{Z})$ are an open algebraic subset of $G=SL_n(\mathbb{Z}/p\mathbb{Z})$.   We apply a theorem of Borel \cite[Theorem B]{Borel} that  states the word map  from $G^d$ to $G$ is dominant in order to show that the set of bad generators are inside a proper sub-variety of $G^d$. We exclude a  proper sub-variety of  generators with such properties  inside $G^d$ that has probability  $O(p^{-1})$. We denote the set of remaining points in  $G^d$ by $\mathcal{F}$. Then we show that for every  $ 1\leq m \leq A$, the variance of  the $m$th moment of monochromatic eigenvalues associated to the remaining generators  is less than $O(p^{-2})$. Recall that 
   \begin{equation*}
 \text{Var}(m):=\frac{1}{|\mathcal{F}|}\mathlarger{\sum_{S\in \mathcal{F}}} (\mu_{\pi,S}(x^m)-N(d,m))^2.
 \end{equation*}
   It follows from the same lines as in the proof of Theorem~\ref{ddd} ; see inequality {\ref{maineq}} that
   \begin{equation}\label{lastone}\text{Var}(m) \leq \frac{2(2d)^{2m}|G|}{\min_{g\in G}(|[g]|) d^2(\pi)},\end{equation}
where  $|[g]|$ denote the size of the conjugacy class of $g\in SL_n(\mathbb{Z}/p\mathbb{Z})$. Note that we give an upper bound on the variance $\text{Var}(m)$ in terms of a group theoretic inequality for $G$.   If $g$ is a regular element of $SL_n(\mathbb{Z}/p\mathbb{Z})$ then $$p^{n^2-n} \ll |[g]|.$$
Let $\pi$ be a non-trivial  irreducible representation of $SL_n(\mathbb{Z}/p\mathbb{Z})$ and $d(\pi)$ be the dimension of this irreducible representation.  Then 
$$p^{n-1}\leq d(\pi).$$
 Since $G=O(p^{n^2-1})$ then
 $$\frac{|G|}{\min_{g}(|[g]|) d^2(\pi)}= O(p^{-(n-1)}).$$
  Hence, by inequality~\ref{lastone}, we deduce that 
 $$ \text{Var}(m) =O\Big(\frac{(2d)^{2m}}{p^{-(n-1)}} \Big).$$
 Since, $n\geq 3$ then 
$$ \text{Var}(m)=O(p^{-2}). $$
 By Chebyshev's inequality, we deduce that
 $$ |\mu_{\pi,S}(x^m)-f_{2d}(x^m)|\geq p^{-1/2},$$
 with probability less than $O(p^{-1})$. This concludes the proof of  inequality~\ref{mme1} and Theorem~\ref{lasttheorem}.

%
%
%

\end{proof}

\subsection*{Acknowledgments}
We would like to  thank Prof. Sarnak for suggesting this project to us and also his insightful comments on the earlier versions of this work.  The computations in this paper were performed using the \textit{Mathematica} system. This material is based upon work supported by the National Science Foundation under Grant No. DMS-1440140 while the second named author was in residence at the Mathematical Sciences Research Institute in Berkeley, California, during the Spring 2017 semester.
\bibliographystyle{alpha}
\bibliography{fixedgenSL2}

\end{document}